\documentclass[1pt]{amsart}

\usepackage{amsmath,amsthm, amscd, amssymb, amsfonts}

\usepackage{hhline}

\newcommand{\Ind}{\operatorname{Ind}}
\newcommand{\ord}{\operatorname{ord}}
\newcommand{\oper}{\#}

\DeclareMathOperator*{\prode}{\boxtimes}

\newcommand\toba{{\mathfrak B }}
\newcommand\trasp{\pi}
\newcommand{\gr}{\operatorname{gr}}
\newcommand{\trid}{\triangleright}

\newcommand\compvert[2]{\genfrac{}{}{-1pt}{0}{#1}{#2}}

\newcommand\mvert[2]{\begin{tiny}\begin{matrix}#1\vspace{-4pt}\\#2\end{matrix}\end{tiny}}
\DeclareMathOperator*{\Tim}{\times}
\newcommand{\Times}[2]{\sideset{_#1}{_#2}\Tim}

\newcommand{\cx}{{\daleth}}

\newcommand{\lu}{{\nu}}
\newcommand{\ld}{{\ell}}

\newcommand{\abofe}{{A(\lu, \ld)}}

\newcommand{\bbofe}{{B(\lu, \ld)}}

\newcommand{\proy}{\Pi}

\newcommand{\Lc}{{\mathcal L}}
\newcommand{\Y}{{\mathcal Y}}
\newcommand{\R}{{\mathcal R}}
\newcommand{\W}{{\mathcal W}}
\newcommand{\Zc}{{\mathcal Z}}
\newcommand{\daga}{{\dagger}}

\newcommand{\acts}{{\rightharpoondown}}

\newcommand{\ku}{\mathbb C}
\newcommand{\K}{{\mathcal K}}
\newcommand{\Z}{{\mathbb Z}}
\newcommand{\N}{{\mathbb N}}
\newcommand{\I}{{\mathbb I}}

\newcommand{\G}{{\mathbb G}}
\newcommand{\M}{{\mathcal M}}
\newcommand{\Q}{{\mathcal Q}}
\newcommand{\F}{{\mathcal F}}
\newcommand{\C}{{\mathcal C}}
\newcommand{\D}{{\mathcal D}}
\newcommand{\Ec}{{\mathbf E}}
\newcommand{\Ee}{{\mathcal E}}
\newcommand{\Kc}{{\mathbf K}}

\newcommand{\uno}{{\bf 1}}
\newcommand{\uv}{{\bf 1_v}}
\newcommand{\uh}{{\bf 1_h}}
\newcommand{\m}{\mathcal{M}}
\newcommand{\n}{\mathcal{N}}

\newcommand{\Dc}{{\mathbf D}}
\newcommand{\B}{{\mathcal B}}
\newcommand{\T}{{\mathcal T}}
\newcommand{\Hc}{{\mathcal H}}
\newcommand{\Vc}{{\mathcal V}}
\newcommand{\Pc}{{\mathcal P}}
\newcommand{\Oc}{{\mathcal O}}
\newcommand{\ydh}{{}^H_H\mathcal{YD}}

\newcommand{\Ss}{{\mathcal S}}
\newcommand{\Vect}{\operatorname{Vec}}
\newcommand{\End}{\operatorname{End}}
\newcommand{\Aut}{\operatorname{Aut}}
\newcommand{\Int}{\operatorname{Int}}
\newcommand{\Ext}{\operatorname{Ext}}
\newcommand\tr{\operatorname{tr}}
\newcommand\Rep{\operatorname{Rep}}

\newcommand\card{\operatorname{card}}
\newcommand{\unosigma}{{\bf 1}^{\sigma}}
\newcommand\tot{\operatorname{tot}}
\newcommand\sgn{\operatorname{sgn}}
\newcommand\ad{\operatorname{ad}}
\newcommand\Hom{\operatorname{Hom}}
\newcommand\opext{\operatorname{Opext}}
\newcommand\Tot{\operatorname{Tot}}
\newcommand\Map{\operatorname{Map}}
\newcommand{\fde}{{\triangleright}}
\newcommand{\fiz}{{\triangleleft}}
\newcommand{\wC}{\widehat{C}}
\newcommand{\la}{\langle}
\newcommand{\ra}{\rangle}
\newcommand{\Comod}{\mbox{\rm Comod\,}}
\newcommand{\Mod}{\mbox{\rm Mod\,}}
\newcommand{\Sg}{{\mathfrak S}}
\newcommand{\funcion}{\mathfrak p}
\newcommand{\tauo}{{\widehat\tau}}
\newcommand{\prin}{t}
\newcommand{\fin}{b}
\newcommand{\pri}{r}
\newcommand{\fine}{l}
\newcommand\rh{\sim_{H}}
\newcommand\rv{\sim_{V}}
\newcommand\rd{\sim_{D}}

\newcommand\V{\operatorname{Vec}}

\theoremstyle{plain}

\renewcommand{\baselinestretch}{1.2}

\title{The character tables of  centralizers in  Weyl Groups of $E_6$, $E_7$, $F_4$, $G_2$}
\author{ \small Shouchuan Zhang,    \ \ Peng Wang, \ \ Jing Cheng,\ \ Hui Yang
}
\address{ Mathematics Department of Hunan University,\newline \indent Changsha China,
410082, E-mail: z9491@yahoo.com.cn }
\date{}

\begin{document}
\newtheorem{Proposition}{\quad Proposition}[section]
\newtheorem{Theorem}{\quad Theorem}
\newtheorem{Definition}[Proposition]{\quad Definition}
\newtheorem{Corollary}[Proposition]{\quad Corollary}
\newtheorem{Lemma}[Proposition]{\quad Lemma}
\newtheorem{Example}[Proposition]{\quad Example}

\maketitle \addtocounter{section}{-1}

\numberwithin{equation}{section}

\date{}

\begin {abstract}  To classify the finite dimensional pointed Hopf
algebras with Weyl group $G$ of $E_6$, $E_7$, $F_4$, $G_2$, we
obtain  the representatives of conjugacy classes of $G$ and  all
character tables of centralizers of these representatives by means
of software GAP.


keywords: GAP, Hopf algebra, Weyl group, character.
\end {abstract}
\section{Introduction}\label {s0}

This article is to contribute to the classification of
finite-dimensional complex pointed Hopf algebras  with Weyl groups
of $E_6$, $E_7$, $F_4$, $G_2$.
 Many papers are about the classification of finite dimensional
pointed Hopf algebras, for example,  \cite{AS98b, AS02, AS00, AS05,
He06, AHS08, AG03, AFZ08, AZ07,Gr00,  Fa07, AF06, AF07, ZZC04, ZCZ08}.

In these research  ones need  the centralizers and character tables
of groups. In this paper we obtain   the representatives of
conjugacy classes of Weyl groups of $E_6$, $E_7$, $F_4$, $G_2$  and all
character tables of centralizers of these representatives by means
of software GAP.

By the Cartan-Killing classification of simple Lie algebras over complex field the Weyl groups
to
be considered are $W(A_l), (l\geq 1); $
 $W(B_l), (l \geq 2);$
 $W(C_l), (l \geq 2); $
 $W(D_l), (l\geq 4); $
 $W(E_l), (8 \geq l \geq 6); $
 $W(F_l), ( l=4 );$
$W(G_l), (l=2)$.

It is otherwise desirable to do this in view of the importance of Weyl groups in the theories of
Lie groups, Lie algebras and algebraic groups. For example, the irreducible representations of
Weyl groups were obtained by Frobenius, Schur and Young. The conjugace classes of  $W(F_4)$
were obtained by Wall
 \cite {Wa63} and its character tables were obtained by Kondo \cite {Ko65}.
The conjugace classes  and  character tables of $W(E_6),$ $W(E_7)$ and $W(E_8)$ were obtained
by  Frame \cite {Fr51}.
Carter gave a unified description of the conjugace classes of
 Weyl groups of simple Lie algebras  \cite  {Ca72}.

\section {Program}

 By using the following
program in GAP,  we  obtain  the representatives of conjugacy
classes of Weyl groups of $E_6$  and all character tables of
centralizers of these representatives.

gap$>$ L:=SimpleLieAlgebra("E",6,Rationals);;

gap$>$ R:=RootSystem(L);;

 gap$>$ W:=WeylGroup(R);Display(Order(W));

gap $>$ ccl:=ConjugacyClasses(W);;

gap$>$ q:=NrConjugacyClasses(W);; Display (q);

gap$>$ for i in [1..q] do

$>$ r:=Order(Representative(ccl[i]));Display(r);;

$>$ od; gap

$>$ s1:=Representative(ccl[1]);;cen1:=Centralizer(W,s1);;

gap$>$ cl1:=ConjugacyClasses(cen1);

gap$>$ s1:=Representative(ccl[2]);;cen1:=Centralizer(W,s1);;

 gap$>$
cl2:=ConjugacyClasses(cen1);

gap$>$ s1:=Representative(ccl[3]);;cen1:=Centralizer(W,s1);;

gap$>$ cl3:=ConjugacyClasses(cen1);

 gap$>$
s1:=Representative(ccl[4]);;cen1:=Centralizer(W,s1);;

 gap$>$
cl4:=ConjugacyClasses(cen1);

 gap$>$
s1:=Representative(ccl[5]);;cen1:=Centralizer(W,s1);;

gap$>$ cl5:=ConjugacyClasses(cen1);

 gap$>$
s1:=Representative(ccl[6]);;cen1:=Centralizer(W,s1);;

gap$>$ cl6:=ConjugacyClasses(cen1);

gap$>$ s1:=Representative(ccl[7]);;cen1:=Centralizer(W,s1);;

gap$>$ cl7:=ConjugacyClasses(cen1);

gap$>$ s1:=Representative(ccl[8]);;cen1:=Centralizer(W,s1);;

gap$>$ cl8:=ConjugacyClasses(cen1);

gap$>$ s1:=Representative(ccl[9]);;cen1:=Centralizer(W,s1);;

gap$>$ cl9:=ConjugacyClasses(cen1);

gap$>$ s1:=Representative(ccl[10]);;cen1:=Centralizer(W,s1);;

 gap$>$
cl10:=ConjugacyClasses(cen1);

gap$>$ s1:=Representative(ccl[11]);;cen1:=Centralizer(W,s1);;

gap$>$ cl11:=ConjugacyClasses(cen1);

gap$>$ s1:=Representative(ccl[12]);;cen1:=Centralizer(W,s1);;

gap$>$ cl2:=ConjugacyClasses(cen1);

gap$>$ s1:=Representative(ccl[13]);;cen1:=Centralizer(W,s1);;

 gap$>$
cl13:=ConjugacyClasses(cen1);

 gap$>$
s1:=Representative(ccl[14]);;cen1:=Centralizer(W,s1);;

gap$>$ cl14:=ConjugacyClasses(cen1);

gap$>$ s1:=Representative(ccl[15]);;cen1:=Centralizer(W,s1);;

gap$>$ cl15:=ConjugacyClasses(cen1);

gap$>$ s1:=Representative(ccl[16]);;cen1:=Centralizer(W,s1);;

gap$>$ cl16:=ConjugacyClasses(cen1);

gap$>$ s1:=Representative(ccl[17]);;cen1:=Centralizer(W,s1);;

gap$>$ cl17:=ConjugacyClasses(cen1);

gap$>$ s1:=Representative(ccl[18]);;cen1:=Centralizer(W,s1);;

gap$>$ cl18:=ConjugacyClasses(cen1);

gap$>$ s1:=Representative(ccl[19]);;cen1:=Centralizer(W,s1);;

gap$>$ cl19:=ConjugacyClasses(cen1);

gap$>$ s1:=Representative(ccl[20]);;cen1:=Centralizer(W,s1);;

gap$>$ cl20:=ConjugacyClasses(cen1);

gap$>$ s1:=Representative(ccl[21]);;cen1:=Centralizer(W,s1);;

 gap$>$
cl21:=ConjugacyClasses(cen1);

gap$>$ s1:=Representative(ccl[22]);;cen1:=Centralizer(W,s1);;

 gap$>$
cl22:=ConjugacyClasses(cen1);

gap$>$ s1:=Representative(ccl[23]);;cen1:=Centralizer(W,s1);;

gap$>$ cl23:=ConjugacyClasses(cen1);

gap$>$ s1:=Representative(ccl[24]);;cen1:=Centralizer(W,s1);;

$>$ cl24:=ConjugacyClasses(cen1);

gap$>$ s1:=Representative(ccl[25]);;cen1:=Centralizer(W,s1);

gap$>$ cl25:=ConjugacyClasses(cen1);

gap$>$ for i in [1..q] do

$>$ s:=Representative(ccl[i]);;cen:=Centralizer(W,s);;

$>$ char:=CharacterTable(cen);;Display (cen);Display(char);

 $>$ od;

The programs for Weyl groups of $E_7$, $E_8$, $F_4$ and $ G_2$ are
similar.

It is similar to $E_6$ about $E_7$, $F_4, G_2$.

\section {$E_6$}

In this section $G$ denotes the Weyl group $W(E_6)$ of $E_6$. It is
clear that $G$  is a  sub-group of general linear group ${\rm GL}
(6, {\mathbb C} )$, \noindent \noindent where ${\mathbb C}$ denotes
the complex field.

The generators of $G$ are:\\
$\left(
\begin{array}{cccccc}-1& 0& 1& 0& 0& 0 \\  0& 1& 0& 0& 0& 0\\
0& 0& 1& 0& 0& 0 \\ 0& 0& 0& 1& 0& 0\\ 0& 0& 0& 0& 1& 0 \\  0& 0&
0&0& 0& 1 \end{array}\right)$,
  $\left(\begin{array}{cccccc}1& 0& 0& 0& 0& 0\\ 0& -1& 0& 1& 0& 0\\ 0& 0& 1& 0& 0& 0\\ 0& 0& 0& 1& 0& 0\\ 0& 0& 0& 0& 1& 0 \\ 0& 0& 0& 0& 0& 1\end{array}\right)$,
   $\left(\begin{array}{cccccc} 1& 0& 0& 0& 0& 0\\ 0& 1& 0& 0& 0& 0\\ 1& 0& -1& 1& 0& 0 \\ 0& 0& 0& 1& 0&
   0\\ 0& 0& 0& 1& 0 \\ 0& 0& 0& 0& 0& 1 \end{array}\right)$,\\
   $\left(\begin{array}{cccccc}1& 0& 0& 0& 0& 0\\ 0& 1& 0& 0& 0& 0\\ 0& 0& 1& 0& 0& 0 \\ 0& 1& 1& -1& 1&\\0& 0& 0& 0& 1& 0 \\ 0& 0& 0& 0& 0& 1 \end{array}\right)$,
  $\left(\begin{array}{cccccc}1& 0& 0& 0& 0& 0\\ 0& 1& 0& 0& 0& 0\\ 0& 0& 1& 0& 0& 0\\ 0& 0& 0& 1& 0& 0\\ 0& 0& 0& 1& -1& 1\\ 0& 0&0& 0& 0& 1\end{array}\right)$,
 $\left(\begin{array}{cccccc} 1&0& 0& 0& 0& 0\\ 0& 1& 0& 0& 0& 0\\ 0& 0& 1& 0& 0&0\\ 0& 0& 0& 1& 0& 0\\ 0& 0& 0& 0& 1& 0\\ 0& 0& 0& 0& 1& -1
 \end{array}\right)$.

The representatives of conjugacy classes of $G$ are:\\
$s_1:= \left(\begin{array}{cccccc} 1& 0& 0& 0& 0& 0 \\ 0& 1& 0& 0& 0& 0
\\
0& 0& 1& 0& 0& 0 \\ 0& 0& 0& 1& 0& 0 \\
     0& 0& 0& 0& 1& 0 \\ 0& 0& 0& 0& 0& 1  \end{array}\right)$,
  $s_2:= \left(\begin{array}{cccccc}-1& 1& 1& -1& 0& 0 \\ -1& 1& 1& -1& 0& 1 \\ -2& 1& 2& -1& 0& 0 \\ -3& 1& 2& -1& 0& 1 \\ -2& 1& 1& -1& 1& 0 \\ -1& 0& 1& -1& 1& 0 \end{array}\right)$,
  $s_3:= \left(\begin{array}{cccccc}1& 0& 0& 0& 0& 0 \\ 0& 0& 1& -1& 1& 0 \\ 0& 0& 1& 0& 0& 0 \\ 0& -1& 1& 0& 1& 0 \\0& 0& 0& 0& 1& 0 \\ 0& 0& 0& 0& 1& -1
  \end{array}\right)$,\\
  $s_4:= \left(\begin{array}{cccccc}-1& 1& 1& -1& 0& 0 \\ 0& 1& 1& -1& -1& 1 \\ 0& 2& 1& -2& 0& 0 \\ 0& 2& 1& -2& -1& 1 \\0& 1& 1& -1& -1& 0 \\ 0& 0& 1& -1& 0& 0 \ \end{array}\right)$,
  $s_5:= \left(\begin{array}{cccccc}0& 0& 1& 0& -1& 0 \\ 0& 1& 1& -1& -1& 1 \\ 1& 1& 1& -1& -1& 0 \\ 1& 1& 1& -1& -2& 1 \\ 0& 1& 1& -1& -1& 0 \\ 0& 0& 1& -1& 0& 0  \end{array}\right)$,
  $s_6:= \left(\begin{array}{cccccc}-1& 1& 1& -1& 0& 0 \\ 0& 1& 1& -2& 1& 0 \\ 0& 2& 1& -2& 0& 0 \\ 0& 2& 1& -3& 1& 0 \\ 0& 2& 0& -2& 1& 0 \\ 0& 1& 0& -1& 1& -1
  \end{array}\right)$,\\
 $s_7:= \left(\begin{array}{cccccc}1& 0& -1& 0& 1& -1 \\ 1& 1& -1& -1& 1& 0 \\ 1& 0& -1& -1& 2& -1\\ 2& 1& -1& -2& 2& -1 \\ 1& 1& 0& -2& 2& -1 \\ 0& 1& 0& -1& 1& -1  \end{array}\right)$,
  $s_8:= \left(\begin{array}{cccccc} 1& -1& 0& 0& 0& 0 \\ 1& 0& 0& 0& -1& 0 \\ 2& 0& 0& -1& 0& 0 \\ 2& 0& 1& -1& -1& 0\\2& 0& 0& 0& -1& 0 \\ 1& 0& 0& 0& -1& 1  \end{array}\right)$,
  $s_9:= \left(\begin{array}{cccccc}-1& 0& 0& 0& 1& -1 \\ -1& -1& 0& 0& 1& 0 \\ -1& -1& -1& 1& 1& -1\\-1& -1& -1& 0& 2& -1 \\ -1& -1& 0& 0& 1& -1 \\ -1& 0& 0& 0& 0& 0
  \end{array}\right)$,\\
  $s_{10}:= \left(\begin{array}{cccccc} 0& 1& 0& 0& 0& -1 \\ 0& 1& 0& 0& 0& 0\\ -1& 1& 1& 0& 0& -1\\ -1& 1& 0& 1& 0& -1 \\-1& 1& 0& 0& 1& -1 \\ -1& 1& 0& 0& 0& 0\end{array}\right)$,
  $s_{11}:= \left(\begin{array}{cccccc} 0& 1& 0& 0& 0& -1 \\ 0& 1& 0& 0& 0& 0 \\ 0& 1& -1& 1& 0& -1 \\ 0& 2& -1& 0& 1& -1\\ -1& 1& 0& 0& 1& -1 \\ -1& 1& 0& 0& 0& 0  \end{array}\right)$,
  $s_{12}:= \left(\begin{array}{cccccc}1& 0& 0& 0& 0& 0 \\ 0& 1& 0& 0& 0& 0 \\ 1& 1& 0& -1& 1& 0 \\ 0& 1& 1& -1& 1& 0 \\ 0& 0& 0& 0& 1& 0 \\ 0& 0& 0& 0& 0& 1
  \end{array}\right)$,\\
  $s_{13}:= \left(\begin{array}{cccccc}0& 0& 0& 0& 0& -1 \\ 0& -1& 1& 0& 0& -1 \\ 0& -1& 0& 1& 0& -2 \\ 0& -1& 1& 0& 1& -3\\-1& -1& 1& 0& 1& -2 \\ -1& 0& 1& 0& 0& -1  \end{array}\right)$,
 $s_{14}:=\left(\begin{array}{cccccc}1& 0& -1& 0& 0& 1 \\ 1& 0& -2& 1& 0& 0 \\ 2& 0& -2& 0& 1& 0 \\ 2& -1& -3& 1& 1& 0 \\ 1& -1& -2& 1& 1& 0 \\ 1& -1& -1& 1& 0& 0 \end{array}\right)$,
  $s_{15}:=\left(\begin{array}{cccccc}1& 0& 0& 0& 0& 0 \\ 0& 0& 1& 0& 0& -1 \\ 1& 0& 1& 0& -1& 0 \\ 1& -1& 1& 1& -1& -1\\0& -1& 1& 1& -1& 0 \\ 0& 0& 0& 1& -1& 0
  \end{array}\right)$,\\
 $s_{16}:=\left(\begin{array}{cccccc}-1& 1& 1& -1& 0& 0 \\ 0& 1& 1& -2& 1& 0 \\ -1& 1& 2& -2& 0& 1 \\ -1& 1& 2& -3& 1& 1 \\-1& 0& 2& -2& 1& 0 \\ 0& 0& 1& -1& 0& 0  \end{array}\right)$,
  $s_{17}:=\left(\begin{array}{cccccc}-1& 0& 0& 1& -1& 0 \\ -1& -1& 0& 1& 0& -1 \\ -1& 0& 0& 1& -1& -1\\-2& 0& 0& 1& 0& -2 \\ -1& 0& -1& 1& 0& -1 \\ -1& 0& 0& 0& 0& 0  \end{array}\right)$,
  $s_{18}:=\left(\begin{array}{cccccc}1& -1& 0& 0& 0& 0 \\ 1& -1& 0& 0& 0& 1 \\ 1& -2& 0& 1& -1& 1 \\ 2& -2& -1& 1& -1& 2\\1& -2& -1& 1& 0& 1 \\ 1& -1& -1& 1& 0& 0
  \end{array}\right)$,\\
 $s_{19}:=\left(\begin{array}{cccccc}-1& 0& 0& 0& 0& 1 \\ -1& -1& 1& 0& 0& 1 \\ -1& 0& 1& -1& 0& 2\\-2& -1& 2& -1& 0& 2 \\ -1& -1& 1& -1& 1& 1 \\ 0& -1& 0& 0& 0& 1  \end{array}\right)$,
 $s_{20}:=\left(\begin{array}{cccccc}-1& 0& 1& 0& 0& 0 \\ -1& 1& 1& -1& 0& 1 \\ -1& 1& 1& -1& 1& 0 \\ -1& 2& 1& -1& 0& 1\\-1& 1& 0& 0& 0& 1 \\ -1& 1& 0& 0& 0& 0  \end{array}\right)$,
  $s_{21}:=\left(\begin{array}{cccccc}-1& 0& 0& 0& 0& 1 \\ 0& 0& -1& 0& 0& 1 \\ -1& 1& -1& 0& 0& 1\\-1& 1& -1& -1& 1& 1 \\-1& 1& 0& -1& 0& 1 \\ -1& 0& 0& 0& 0& 0
  \end{array}\right)$,\\
  $s_{22}:=\left(\begin{array}{cccccc}1& -1& 0& 0& 0& 0 \\ 1& -1& -1& 1& -1& 1 \\ 1& -2& 0& 1& -1& 1 \\2& -2& -1& 1& -1& 2 \\ 1& -2& -1& 1& 0& 1\\1& -1& -1& 1& 0& 0 \end{array}\right)$,
  $s_{23}:=\left(\begin{array}{cccccc}-1& 1& 1& -1& 0& 0 \\ -1& 0& 1& 0& -1& 0 \\ -1& 1& 1& -1& -1& 1 \\-2& 1& 1& 0& -2& 1 \\ -1& 1& 0& 0& -1& 0 \\ 0& 1& 0& 0& -1& 0   \end{array}\right)$,
 $s_{24}:=\left(\begin{array}{cccccc}1& -1& 0& 0& 0& 0 \\ 1& -1& 0& 0& 0& 1 \\ 2& -1& 0& 0& 0& 0 \\ 2& -2& 0& 1& -1& 1 \\1& -2& 0& 1& 0& 0 \\ 0& -1& 1& 0& 0& 0
 \end{array}\right)$,\\
 $s_{25}:=\left(\begin{array}{cccccc}-1& 0& 0& 0& 0& 1 \\ 0& -1& -1& 1& 0& 1 \\ 0& 0& -1& 0& 0& 2 \\ 0& -1& -1& 1& -1& 3\\0& -1& -1& 1& -1& 2 \\ 0& -1& 0& 0& 0& 1  \end{array}\right)$

We  orderly denote  these representatives  above  by  $s_1, s_2, \cdots, s_{25}$.
 Obviously,  $s_1$ is the unity element
and
$G= G^{s_1}$.
The character table of $G^{s_1} =G$:\\
{ \small
 \noindent \begin{tabular}{c|ccccccccccccccccccccccccc}
    & & & & & & & & & &10 & & & & & & & & & &20 & & & & &                                    \\\hline
  $\chi_  {1}^ {(1)}$   &  1 &1&1& 1&1&1&1&1&1& 1&1&1&1&1&1&1&1& 1&1&1&1& 1&1& 1&1\\

  $\chi _{1}^ {(2)}$   &1&-1&1&-1&1 &-1&1 &-1& -1& 1& -1&1&1& -1& -1&1& -1&-1&1&1&1& 1&1& 1&1\\

  $\chi _{1}^ {(3)}$   &  6& 1&1&-4&.&2&3&.&.&-2 &-2&2&2& -1&2&.&.&-1 &-1& -2&.&-3&1&-1&1\\

  $\chi _{1}^ {(4)}$    &6&-1&1& 4&. &-2&3&.&.&-2&2&2&2&1& -2&.&.& 1 &-1&-2&.&-3&1&-1 &1\\

  $\chi _{1}^ {(5)}$   &10& .&.& .&4&. &-2&.&.&-6&.&2&2&.&.& -2&.& .&2&.&1& 1& -3&-1&.\\

  $\chi _{1}^ {(6)}$   &15& .&.&-5&.& -2&3&.&3&-1& -1 &-1&3&1&-1& -1&1&-1& -1&2&.& 6&2& .&-1\\

  $\chi _{1}^ {(7)}$  &15& .&.&-5&3&1&.& -1 &-1& 7&1&3 &-1 &-2&-3&1&1& .&.&1&.&-3&1&-1&-2\\
\end{tabular} }
{\bf Remark:} The top row of the table above should fill  the representatives $s_1, s_2, \cdots, s_{25}$
of conjugacy classes of $G^{s_1} =G$. Considering the width of the table we omit them.
 we omit the top rows in the following character tables   for the same reason. We write 10, 20 , $\cdots$, to show
 Column 10, Column 20, $\cdots$.


The generators of $G^{s_2}$ are:\\
$\left(\begin{array}{cccccc}-1& 1& 1& -1& 0& 0 \\ -1& 1& 1& -1& 0& 1
\\ -2& 1& 2& -1& 0& 0 \\ -3& 1& 2& -1& 0& 1\\ -2& 1& 1& -1& 1& 0 \\
-1& 0& 1& -1& 1& 0\end{array}\right)$,$\left(\begin{array}{cccccc}1&
0& 0& 0& 0& 0
\\ 1& 0& 1& -1& 0& 0\\ 2& 1& 0& -1& 0& 0 \\ 3& 0& 0& -1& 0& 0\\2& 0& 0& 0& -1& 0 \\ 1& 0& 0& 0& 0& -1\end{array}\right)$,
$\left(\begin{array}{cccccc}-1& 1& 1& -1& 0& 0 \\ 0& 2& 0& -1& 0& 0
\\ 0& 2& 1& -2& 0& 0 \\ 0& 3& 0& -2& 0& 0\\0& 2& 0& -2& 1& 0 \\ 0& 1& 0& -1& 0& 1
\end{array}\right)$,\\
$\left(\begin{array}{cccccc}1& 0& 0& 0& 0& 0 \\ 0& 0& 1& -1& 1& 0
\\1& 1& 0& -1& 1& 0 \\ 1& 0& 0& -1& 2& 0\\0& 0& 0& 0& 1& 0 \\ 0& 0& 0& 0& 1& -1 \end{array}\right)$.

The representatives of conjugacy classes of   $G^{s_2}$ are:\\
$\left(\begin{array}{cccccc}1& 0& 0& 0& 0& 0 \\ 0& 1& 0& 0& 0& 0
\\0& 0& 1& 0& 0& 0 \\ 0& 0& 0& 1& 0& 0\\0& 0& 0& 0& 1& 0 \\ 0& 0& 0& 0& 0& 1 \end{array}\right)$,
 $\left(\begin{array}{cccccc}-1& 1& 1& -1& 0& 0 \\ -1& 0& 2& -1& 0& 0 \\ -1& 1& 2& -1& -1& 0 \\-2& 0& 3& -1& -1& 0 \\ -2& 0& 2& 0& -1& 0 \\ -1& 0& 1& 0& 0& -1\end{array}\right)$,
$\left(\begin{array}{cccccc}-1& 1& 1& -1& 0& 0 \\ -1& 1& 1& -1& 0& 1
\\ -2& 1& 2& -1& 0& 0\\-3& 1& 2& -1& 0& 1 \\ -2& 1& 1& -1& 1& 0 \\ -1& 0& 1& -1& 1& 0
\end{array}\right)$,\\
  $\left(\begin{array}{cccccc} -1& 1& 1& -1& 0& 0 \\ -1& 1& 1& -1& 0& 1 \\ -1& 2& 1& -2& 1& 0 \\-2& 2& 1& -2& 1& 1 \\ -2& 1& 1& -1& 1& 0 \\ -1& 0& 1& -1& 1& 0 \end{array}\right)$,
$\left(\begin{array}{cccccc}-1& 1& 1& -1& 0& 0 \\ -1& 1& 1& -1& 1&
-1 \\ -2& 1& 2& -1& 0& 0\\ -3& 1& 2& -1& 1& -1 \\ -2& 1& 1& -1& 1& 0
\\ -1& 1& 0& 0& 0& 0 \end{array}\right)$,
$\left(\begin{array}{cccccc}-1& 1& 1& -1& 0& 0 \\ -1& 1& 1& -1& 1&
-1 \\ -1& 2& 1& -2& 1& 0\\-2& 2& 1& -2& 2& -1 \\ -2& 1& 1& -1& 1& 0 \\
-1& 1& 0& 0& 0& 0 \end{array}\right)$,\\
 $\left(\begin{array}{cccccc}-1& 1& 1& -1& 0& 0 \\ -1& 1& 1& 0& -1& 0 \\ -2& 1& 2& -1& 0& 0\\-3& 1& 2& 0& -1& 0 \\ -2& 0& 2& 0& -1& 0 \\ -1& 0& 1& 0& -1& 1
 \end{array}\right)$,
  $\left(\begin{array}{cccccc}-1& 1& 1& -1& 0& 0 \\ -1& 1& 1& 0& -1& 0 \\ -1& 1& 2& -1& -1& 0\\-2& 1& 2& 0& -2& 0 \\ -2& 0& 2& 0& -1& 0 \\ -1& 0& 1& 0& -1&
  1\end{array}\right)$,
$\left(\begin{array}{cccccc}-1& 1& 1& -1& 0& 0 \\ 0& 1& 1& -1& -1& 1
\\ 0& 2& 1& -2& 0& 0 \\ 0& 2& 1& -2& -1& 1\\0& 1& 1& -1& -1& 0 \\ 0& 0& 1& -1& 0& 0
\end{array}\right)$,\\
 $\left(\begin{array}{cccccc} -1& 1& 1& -1& 0& 0 \\ 0& 1& 1& -1& 0& -1 \\ 0& 2& 1& -2& 0& 0 \\ 0& 2& 1& -2&
0& -1 \\ 0& 1& 1& -1& -1& 0 \\ 0& 1& 0& 0& -1&
0\end{array}\right)$, $\left(\begin{array}{cccccc}1& 0& 0& 0& 0&
0 \\ 0& 0& 1& -1& 1& 0 \\0& 0& 1& 0& 0& 0 \\ 0& -1& 1& 0& 1& 0
\\0& 0& 0& 0& 1& 0 \\ 0& 0& 0& 0& 1& -1\end{array}\right)$,
 $\left(\begin{array}{cccccc}1& 0& 0& 0& 0& 0 \\ 0& 0& 1& -1& 1& 0 \\ 1& 1& 0& -1& 1& 0 \\ 1& 0& 0& -1& 2& 0
\\ 0& 0& 0& 0& 1& 0 \\ 0& 0& 0& 0& 1& -1\end{array}\right)$,\\
$\left(\begin{array}{cccccc}1& 0& 0& 0& 0& 0 \\ 0& 0& 1& 0& -1& 1\\
0& 0& 1& 0& 0& 0 \\ 0& -1& 1& 1& -1& 1 \\ 0& -1& 1& 1& -1& 0 \\ 0&
-1& 1& 0& 0& 0 \end{array}\right)$,
 $\left(\begin{array}{cccccc}1& 0& 0& 0& 0& 0 \\ 0& 0& 1& 0& -1& 1 \\ 1& 0& 1& 0& -1& 0 \\ 1& -1& 1& 1& -2&
1 \\0& -1& 1& 1& -1& 0 \\ 0& -1& 1& 0& 0& 0\end{array}\right)$,
 $\left(\begin{array}{cccccc}1& 0& 0& 0& 0& 0 \\ 0& 0& 1& 0& 0& -1 \\ 0& 0& 1& 0& 0& 0 \\ 0& -1& 1& 1& 0& -1
\\ 0& -1& 1& 1& -1& 0 \\ 0& 0& 0& 1& -1& 0\end{array}\right)$,\\
 $\left(\begin{array}{cccccc}1& 0& 0& 0& 0& 0 \\ 0& 0& 1& 0& 0& -1 \\ 1& 0& 1& 0& -1& 0 \\ 1& -1& 1& 1& -1&
-1 \\0& -1& 1& 1& -1& 0 \\ 0& 0& 0& 1& -1& 0\end{array}\right)$,
 $\left(\begin{array}{cccccc}-1& 1& 1& -1& 0& 0 \\ -1& 0& 2& -1& 0& 0 \\ -2& 1& 2& -1& 0& 0 \\-3& 0& 3& -1& 0& 0 \\ -2& 0& 2& 0& -1& 0 \\ -1& 0& 1& 0& 0& -1
 \end{array}\right)$,
$\left(\begin{array}{cccccc} 1& 0& 0& 0& 0& 0 \\ 0& 1& 0& 0& 0&
0\\1& 1& 0& -1& 1& 0 \\ 1& 1& -1& 0& 1& 0 \\ 0& 0& 0& 0& 1& 0 \\ 0&
0& 0& 0& 0&
1\end{array}\right)$,\\$\left(\begin{array}{cccccc}1& 0& 0& 0&
0& 0 \\ 1& 0& 1& -1& 0& 0 \\2& 1& 0& -1& 0& 0 \\ 3& 0& 0& -1& 0&
0\\2& 0& 0& 0& -1& 0 \\ 1& 0& 0& 0& 0& -1\end{array}\right)$,
$\left(\begin{array}{cccccc}1& 0& 0& 0& 0& 0 \\ 1& 1& 0& 0& -1& 0\\
2& 1& 0& -1& 0& 0 \\ 3& 1& -1& 0& -1& 0
\\ 2& 0& 0& 0& -1& 0 \\ 1& 0& 0& 0& -1& 1\end{array}\right).$

The character table of $G^{s_2}$:\\
\noindent \begin{tabular}{c|cccccccccccccccccccc}
  & & & & & & & & & & 10& & & & & & & & & & 20\\\hline

  $\chi _{2}^ {(1)}$    &1&1&1&1&1&1&1&1&1&1&1&1&1&1&1&1&1&1&1
&1\\

  $\chi _{2}^ {(2)}$&    1&1&1 &-1&1 &-1 &-1&1&1&1&1 &-1 &-1&1 &-1&1& -1 &-1&1
&1\\

$\chi _{2}^ {(3)}$   & 1& -1&1&1&1&1 &-1 &-1&1&1&1&1 &-1 &-1 &-1 &-1 &-1&1&1
&1\\

  $\chi _{2}^ {(4)}$    &1&1 &-1&1 &-1&1 &-1&1 &-1 &-1&1 &-1&1 &-1&1 &-1 &-1 &-1&1
&1\\

$\chi _{2}^ {(5)}$&1& -1 &-1 &-1 &-1 &-1 &-1 &-1 &-1 &-1&1&1&1&1&1&1 &-1&1&1
&1\\

$\chi _{2}^ {(6)}$&1 &-1 &-1&1& -1&1&1 &-1 &-1 &-1&1 &-1& -1&1 &-1&1&1& -1&1
&1\\

  $\chi _{2}^ {(7)}$&1&1 &-1 &-1 &-1 &-1&1&1 &-1& -1&1&1 &-1 &-1 &-1 &-1&1&1&1
&1\\

$\chi _{2}^ {(8)}$&1& -1&1& -1&1& -1&1 &-1&1&1&1 &-1&1 &-1&1 &-1&1& -1&1
&1\\

$\chi _{2}^ {(9)}$&1&1&A&A &-A &-A &-1 &-1&A &-A &-1 &-1 &-A &-A&A&A&1&1 &-1
&1\\

$\chi _{2}^ {({10})}$& 1&1 &-A &-A&A&A &-1 &-1 &-A&A &-1 &-1&A&A &-A &-A&1&1& -1
&1\\

$\chi _{2}^ {()}{11}$& 1 &-1&A &-A &-A&A &-1&1&A &-A &-1&1 &-A&A&A &-A&1 &-1 &-1
&1\\

$\chi _{2}^ {({12})}$& 1 &-1 &-A&A&A &-A &-1&1 &-A&A& -1&1&A &-A &-A&A&1 &-1 &-1
&1\\

$\chi _{2}^ {({13})}$& 1& -1 &-A &-A&A&A&1&1 &-A&A &-1 &-1 &-A &-A&A&A &-1&1& -1
&1\\

$\chi _{2}^ {({14})}$& 1 &-1&A&A& -A& -A&1&1&A &-A &-1 &-1&A&A &-A &-A &-1&1& -1
&1\\

$\chi _{2}^ {({15})}$& 1&1 &-A&A&A &-A&1 &-1 &-A&A &-1&1 &-A&A&A &-A &-1 &-1& -1
&1\\

$\chi _{2}^ {({16})}$& 1&1&A& -A &-A&A&1 &-1&A &-A &-1&1&A& -A &-A&A &-1 &-1& -1
&1\\

$\chi _{2}^ {({17})}$& 2&.& -2&.& -2&.&.&.&2&2&2&.&.&.&.&.&.&.& -2
&-2\\

$\chi _{2}^ {({18})}$& 2&.&2&.&2&.&.&. &-2 &-2&2&.&.&.&.&.&.&.& -2
&-2\\

$\chi _{2}^ {({19})}$& 2&.&B&.& -B&.&.&.& -B&B &-2&.&.&.&.&.&.&.& 2
&-2\\

$\chi _{2}^ {({20})}$& 2&. &-B&.&B&.&.&.&B &-B &-2&.&.&.&.&.&.&.& 2
&-2\\
\end{tabular}\\
\noindent \noindent where A = E(4)
  = ER(-1) = i,
B = 2*E(4)
  = 2*ER(-1) = 2i.

The generators of $G^{s_3}$ are:\\
$\left(\begin{array}{cccccc}1& 0& 0& 0& 0& 0 \\ 0& 0& 1& -1& 1&
0\\0& 0& 1& 0& 0& 0 \\ 0& -1& 1& 0& 1& 0\\0& 0& 0& 0& 1& 0 \\ 0& 0&
0& 0& 1& -1\end{array}\right)$, $\left(\begin{array}{cccccc}-1& 0&
1& 0& 0& 0 \\ 0& 1& 0& 0& 0& 0\\ 0& 0& 1& 0& 0& 0 \\ 0& 0& 0& 1& 0&
0
\\0& 0& 0& 0& 1& 0 \\ 0& 0& 0& 0& 0& 1\end{array}\right)$,
 $\left(\begin{array}{cccccc}1& 1& 0& -1& 0& 0 \\ 0& 2& 0& -1& 0& 0 \\ 0& 2& 1& -2& 0& 0 \\ 0& 3& 0& -2& 0& 0\\0& 2& 0& -2& 1& 0 \\ 0& 1& 0& -1& 0& 1
 \end{array}\right)$,\\
$\left(\begin{array}{cccccc}1& 0& 0& 0& 0& 0 \\ 1& 0& 1& -1& 0&
0\\2& 1& 0& -1& 0& 0 \\ 3& 0& 0& -1& 0& 0 \\2& 0& 0& 0& -1& 0 \\ 1&
0& 0& 0& 0& -1 \end{array}\right)$, $\left(\begin{array}{cccccc} 1&
0& 0& 0& 0& 0 \\ 0& 1& 0& 0& 0& 0\\ 0& 0& 1& 0& 0& 0 \\ 0& 0& 0& 1&
0& 0\\ 0& 0& 0& 0& 1& 0 \\ 0& 0& 0& 0& 1& -1\end{array}\right)$,
$\left(\begin{array}{cccccc}-1& 1& 1& -1& 0& 0 \\ 0& 1& 1& -1& 0& -1
\\ 0& 2& 1&
 -2& 0& 0 \\ 0& 2& 1& -2& 0& -1\\0& 1& 1& -1& -1& 0 \\ 0& 0& 1& -1& 0& 0
 \end{array}\right)$,\\
$\left(\begin{array}{cccccc}0& 0& 1& 0& -1& 0 \\ 1& 1& 0& 0& -1& 0
\\ 1& 1& 1& -1& -1& 0 \\ 2& 1& 0& 0& -2& 0\\2& 0& 0& 0& -1& 0 \\ 1& 0& 0& 0& -1& 1\end{array}\right)$.

The representatives   of conjugacy classes of   $G^{s_3}$

$\left(\begin{array}{cccccc}1& 0& 0& 0& 0& 0 \\ 0& 1& 0& 0& 0& 0\\
0& 0& 1& 0& 0& 0 \\ 0& 0& 0& 1& 0& 0 \\0& 0& 0& 0& 1& 0 \\ 0& 0& 0&
0& 0& 1\end{array}\right)$,
 $\left(\begin{array}{cccccc}-1& 0& 1& 0& 0& 0 \\ -1& 0& 2& -1& 0& 0 \\ -2& 1& 2& -1& 0& 0 \\ -3& 0& 3& -1&
0& 0\\-2& 0& 2& 0& -1& 0 \\ -1& 0& 1& 0& 0&
-1\end{array}\right)$, $\left(\begin{array}{cccccc}-1& 0& 1& 0&
0& 0 \\ -1& 0& 2& -1& 0& 0
\\ -1& 0& 2& 0& -1& 0\\ -2& -1& 3& 0& -1& 0 \\ -2& 0& 2& 0& -1& 0 \\ -1& 0& 1& 0& -1&
1\end{array}\right)$,\\
$\left(\begin{array}{cccccc}-1& 0& 1& 0& 0& 0 \\ -1& 0& 2& -1& 0& 0
\\ -1& 0& 2& 0& -1& 0 \\ -2& -1& 3& 0& -1& 0 \\ -2& 0& 2& 0& -1& 0 \\ -1& 0& 1& 0& 0& -1\end{array}\right)$,
$\left(\begin{array}{cccccc}-1& 0& 1& 0& 0& 0 \\ -1& 1& 1& -1& 0& 1
\\ -2& 1& 2& -1& 0& 0 \\ -3& 1& 2& -1& 0& 1\\-2& 1& 1& -1& 1& 0 \\ -1& 0& 1& -1& 1& 0\end{array}\right)$,
$\left(\begin{array}{cccccc}-1& 0& 1& 0& 0& 0 \\ -1& 1& 1& -1& 0& 1
\\ -2& 1& 2& -1& 0& 0 \\ -3& 1& 2& -1& 0& 1\\-2& 1& 1& -1& 1& 0 \\ -1& 1& 0& 0& 0&
0\end{array}\right)$,\\
$\left(\begin{array}{cccccc}-1& 0& 1& 0& 0& 0 \\ -1& 1& 1& -1& 0& 1
\\ -1& 1& 1& -1& 1& 0 \\ -2& 1& 1& -1& 1& 1\\-2& 1& 1& -1& 1& 0 \\ -1& 0& 1& -1& 1& 0\end{array}\right)$,
$\left(\begin{array}{cccccc}-1& 0& 1& 0& 0& 0 \\ -1& 1& 1& -1& 0& 1
\\ -1& 1& 1& -1& 1& 0 \\ -2& 1& 1& -1& 1& 1\\-2& 1& 1& -1& 1& 0 \\ -1& 1& 0& 0& 0& 0 \end{array}\right)$,
$\left(\begin{array}{cccccc}-1& 0& 1& 0& 0& 0 \\ -1& 1& 1& 0& -1& 0
\\ -2& 1& 2& -1& 0& 0 \\ -3& 1& 2& 0& -1& 0\\-2& 0& 2& 0& -1& 0 \\ -1& 0& 1& 0& -1&
1\end{array}\right)$,\\
$\left(\begin{array}{cccccc}-1& 0& 1& 0& 0& 0 \\ -1& 1& 1& 0& -1& 0
\\ -1& 0& 2& 0& -1& 0 \\ -2& 0& 2& 1& -2& 0\\-2& 0& 2& 0& -1& 0 \\ -1& 0& 1& 0& -1& 1 \end{array}\right)$,
$\left(\begin{array}{cccccc}-1& 0& 1& 0& 0& 0 \\ 0& 0& 1& -1& 1& 0
\\ 0& 0& 1& 0& 0& 0 \\ 0& -1& 1& 0& 1&
 0\\0& 0& 0& 0& 1& 0 \\ 0& 0& 0& 0& 0& 1\end{array}\right)$,
$\left(\begin{array}{cccccc}-1& 0& 1& 0& 0& 0 \\ 0& 0& 1& -1& 1& 0
\\ 0& 0& 1& 0& 0& 0 \\ 0& -1& 1& 0& 1& 0 \\0& 0& 0& 0& 1& 0 \\ 0& 0& 0& 0& 1& -1
\end{array}\right)$,\\
$\left(\begin{array}{cccccc}-1& 0& 1& 0& 0& 0 \\ 0& 0& 1& 0& -1& 1
\\ 0& 0& 1& 0& 0& 0 \\ 0& -1& 1& 1& -1& 1\\0& -1& 1& 1& -1& 0 \\ 0& -1& 1& 0& 0& 0 \end{array}\right)$,
$\left(\begin{array}{cccccc}-1& 0& 1& 0& 0& 0 \\ 0& 0& 1& 0& -1& 1
\\ 0& 0& 1& 0& 0& 0 \\ 0& -1& 1& 1& -1& 1\\0& -1&
 1& 1& -1& 0 \\ 0& 0& 0& 1& -1& 0\end{array}\right)$,
$\left(\begin{array}{cccccc}-1& 0& 1& 0& 0& 0 \\ 0& 1& 0& 0& 0&
0\\0& 0& 1& 0& 0& 0 \\ 0& 0& 0& 1& 0& 0 \\ 0& 0& 0& 0& 1& 0 \\ 0& 0&
0& 0& 0& 1\end{array}\right)$, \\$\left(\begin{array}{cccccc}-1&
1& 1& -1& 0& 0 \\ -1& 0& 2& -1& 0& 0
\\ -2& 1& 2& -1& 0& 0\\-3& 0& 3& -1& 0& 0 \\ -2& 0& 2& 0& -1& 0 \\ -1& 0& 1& 0& -1& 1 \end{array}\right)$,
$\left(\begin{array}{cccccc}-1& 1& 1& -1& 0& 0 \\ -1& 0& 2& -1& 0& 0
\\ -2& 1& 2& -1& 0& 0 \\-3& 0& 3& -1& 0& 0 \\ -2& 0& 2& 0& -1& 0 \\ -1& 0& 1& 0& 0& -1\end{array}\right)$,
$\left(\begin{array}{cccccc}-1& 1& 1& -1& 0& 0 \\ -1& 1& 1& -1& 0& 1
\\ -2& 1& 2& -1& 0& 0\\-3& 1& 2& -1& 0& 1 \\ -2& 1& 1& -1& 1&
 0 \\ -1& 0& 1& -1& 1& 0\end{array}\right)$,\\
$\left(\begin{array}{cccccc} -1& 1& 1& -1& 0& 0 \\ -1& 1& 1& -1& 0&
1 \\ -2& 1& 2& -1& 0& 0\\-3& 1& 2& -1& 0& 1 \\ -2& 1& 1& -1& 1& 0 \\
-1& 1& 0& 0& 0& 0\end{array}\right)$,
$\left(\begin{array}{cccccc}-1& 1& 1& -1& 0& 0 \\ -1& 1& 1& 0& -1& 0
\\ -2& 1& 2& -1& 0& 0\\-3& 1& 2& 0& -1& 0 \\ -2& 0& 2& 0& -1& 0 \\ -1& 0& 1& 0& -1& 1\end{array}\right)$,
$\left(\begin{array}{cccccc}-1& 1& 1& -1& 0& 0 \\ 0& 1& 1& -1& -1& 1
\\ 0& 2& 1& -2& 0& 0 \\ 0& 2& 1& -2& -1& 1\\0& 1& 1& -1& -1& 0 \\ 0& 0& 1& -1& 0&
0\end{array}\right)$,\\
$\left(\begin{array}{cccccc}-1& 1& 1& -1& 0& 0 \\ 0& 1& 1& -1& -1& 1
\\ 0& 2& 1& -2& 0& 0 \\ 0& 2& 1& -2& -1& 1\\0& 1& 1& -1& -1& 0 \\ 0& 1& 0& 0& -1& 0\end{array}\right)$,
$\left(\begin{array}{cccccc}1& 0& 0& 0& 0& 0 \\ 0& 0& 1& -1& 1&
0\\0& 0& 1& 0& 0& 0 \\ 0& -1& 1& 0& 1& 0 \\ 0& 0& 0& 0& 1& 0 \\ 0&
0& 0& 0& 0& 1\end{array}\right)$, $\left(\begin{array}{cccccc}1&
0& 0& 0& 0& 0 \\ 0& 0& 1& -1& 1& 0 \\0& 0& 1& 0& 0& 0 \\ 0& -1& 1&
0& 1& 0\\0& 0& 0& 0& 1& 0 \\ 0& 0& 0& 0& 1&
-1\end{array}\right)$,\\ $\left(\begin{array}{cccccc}-1& 0& 1&
0& 0& 0 \\ -1& 0& 2& -1& 0& 0
\\ -2& 1& 2& -1& 0& 0 \\ -3& 0& 3& -1& 0& 0\\-2& 0& 2& 0& -1& 0 \\ -1& 0& 1& 0& -1& 1\end{array}\right).$

The character table of $G^{s_3}$:\\
\noindent
\begin{tabular}{c|ccccccccccccccccccccccccc}
  & & & & & & & & & &10 & & & & & & & & & &20 & & & & &\\\hline
  $\chi _{3}^ {(1)}$      &1  &1  &1  &1  &1  &1   &1  &1  &1  &1  &1  &1  &1  &1  &1  &1  &1  &1  &1  &1  &1  &1&1  &1
  &1\\

  $\chi _{3}^ {(2)}$      &1&-1&-1 &1&-1 &1  &1&-1&-1 &1 &1&-1&-1 &1&-1&-1 &1 &1&-1 &1 &1&-1&-1 &1
 &1\\

  $\chi _{3}^ {(3)}$     &1&-1 &1 &1 &1 &1 &-1&-1&-1 &1&-1&-1 &1 &1&-1 &1 &1&-1&-1 &1&-1&-1 &1 &1
 &-1\\

  $\chi _{3}^ {(4)}$     &1&-1&-1 &1 &1&-1 &-1 &1&-1 &1 &1&-1 &1&-1&-1&-1 &1&-1 &1 &1&-1 &1&-1 &1
 &1\\

  $\chi _{3}^ {(5)}$     &1&-1 &1 &1&-1&-1  &1 &1&-1 &1&-1&-1&-1&-1&-1 &1 &1 &1 &1 &1 &1 &1&1 &1
 &-1\\

  $\chi _{3}^ {(6)}$     &1 &1&-1 &1 &1&-1  &1&-1 &1 &1&-1 &1 &1&-1 &1&-1 &1 &1&-1 &1 &1&-1&-1 &1
 &-1\\

  $\chi _{3}^ {(7)}$     &1 &1 &1 &1&-1&-1 &-1&-1 &1 &1 &1 &1&-1&-1 &1 &1 &1&-1&-1 &1&-1&-1&1 &1
 &1\\
\end{tabular}


The generators of $G^{s_4}$ are:\\

$\left(\begin{array}{cccccc}-1& 1& 1& -1& 0& 0 \\ 0& 1& 1& -1& -1& 1
\\ 0& 2& 1& -2& 0& 0 \\ 0& 2& 1& -2& -1& 1 \\ 0& 1& 1& -1& -1&
0 \\ 0& 0& 1& -1& 0& 0 \end{array}\right)$,
$\left(\begin{array}{cccccc} -1& 0& 1& 0& 0& 0 \\ 0& 1& 0& 0& 0& 0
\\ 0& 0& 1&
 0& 0& 0 \\ 0& 0& 0& 1& 0& 0\\0& 0& 0& 0& 1& 0 \\ 0& 0& 0& 0& 0& 1 \end{array}\right)$,
$\left(\begin{array}{cccccc}1& 1& 0& -1& 0& 0 \\ 0& 2& 0& -1& 0& 0
\\ 0& 2& 1& -2& 0& 0 \\ 0& 3& 0& -2& 0& 0\\0& 2& 0& -2& 1& 0 \\ 0& 1& 0& -1& 0& 1
\end{array}\right)$,\\
$\left(\begin{array}{cccccc}1& 0& 0& 0& 0& 0 \\ 1& 0& 1& -1& 0&
0\\2& 1& 0& -1& 0& 0 \\ 3& 0& 0& -1& 0& 0\\2& 0& 0& 0& -1& 0 \\ 1&
0& 0& 0& 0& -1 \end{array}\right)$, $\left(\begin{array}{cccccc}1&
1& 0& -1& 0& 0 \\ 0& 1& 1& -2& 1& 0
\\ 1& 2& 0& -2& 1& 0 \\ 1& 2& 0& -3& 2& 0\\0& 2& 0& -2& 1& 0 \\ 0& 1& 0& -1& 1& -1 \end{array}\right).$

The representatives   of conjugacy classes of   $G^{s_4}$ are:\\
$\left(\begin{array}{cccccc}1& 0& 0& 0& 0& 0 \\ 0& 1& 0& 0& 0& 0
\\0& 0& 1& 0& 0& 0 \\ 0& 0& 0& 1& 0& 0\\0& 0& 0& 0& 1& 0 \\ 0& 0& 0& 0& 0& 1\end{array}\right)$,
$\left(\begin{array}{cccccc}-1& 0& 1& 0& 0& 0 \\ -1& 0& 2& -1& 0& 0
\\ -1& 0& 2& 0& -1& 0\\-2& -1& 3& 0& -1& 0 \\ -2& 0& 2& 0& -1& 0 \\ -1& 0& 1& 0& 0& -1 \end{array}\right)$,
$\left(\begin{array}{cccccc} -1& 0& 1& 0& 0& 0 \\ -1& 1& 1& -1& 0& 1
\\ -2& 1& 2& -1& 0& 0 \\ -3& 1& 2& -1& 0& 1\\-2& 1& 1& -1& 1& 0 \\ -1& 0& 1& -1& 1&
0\end{array}\right)$,\\
$\left(\begin{array}{cccccc} -1& 0& 1& 0& 0& 0 \\ -1& 1& 1& -1& 0& 1
\\ -1& 1& 1& -1& 1& 0 \\ -2& 1& 1& -1& 1& 1\\-2& 1& 1& -1& 1& 0 \\ -1& 0& 1& -1& 1& 0\end{array}\right)$,
$\left(\begin{array}{cccccc}-1& 0& 1& 0& 0& 0 \\ -1& 1& 1& -1& 1& -1
\\ -2& 1& 2& -1& 0& 0 \\-3& 1& 2& -1& 1& -1 \\ -2& 1& 1& -1& 1& 0 \\ -1& 1& 0& 0& 0&
0\end{array}\right)$,$\left(\begin{array}{cccccc} -1& 0& 1& 0&
0& 0 \\ -1& 1& 1& -1& 1& -1 \\ -1& 1& 1& -1& 1& 0\\-2& 1& 1& -1& 2&
-1
\\ -2& 1& 1& -1& 1& 0 \\ -1& 1& 0&
 0& 0& 0\end{array}\right)$,\\
$\left(\begin{array}{cccccc}-1& 0& 1& 0& 0& 0 \\ -1& 1& 1& 0& -1& 0
\\ -2& 1& 2& -1& 0& 0 \\ -3& 1& 2& 0& -1& 0 \\-2& 0& 2& 0& -1& 0 \\ -1& 0& 1& 0& -1& 1\end{array}\right)$,
$\left(\begin{array}{cccccc}-1& 0& 1& 0& 0& 0 \\ -1& 1& 1& 0& -1& 0
\\ -1& 0& 2& 0& -1& 0 \\ -2& 0& 2& 1& -2& 0\\-2& 0& 2& 0& -1& 0 \\ -1& 0& 1& 0& -1& 1\end{array}\right)$,
$\left(\begin{array}{cccccc}-1& 0& 1& 0& 0& 0 \\ 0& 0& 1& -1& 1& 0
\\ 0& 0& 1& 0& 0& 0 \\ 0& -1& 1& 0& 1& 0 \\0& 0& 0& 0& 1& 0 \\ 0& 0& 0&
 0& 1& -1\end{array}\right)$,\\
$\left(\begin{array}{cccccc}-1& 0& 1& 0& 0& 0 \\ 0& 0& 1& 0& -1& 1
\\ 0& 0& 1& 0& 0& 0 \\ 0& -1& 1& 1& -1& 1 \\0& -1& 1& 1& -1& 0 \\ 0& -1& 1& 0& 0& 0 \end{array}\right)$,
$\left(\begin{array}{cccccc}-1& 0& 1& 0& 0& 0 \\ 0& 0& 1& 0& 0& -1
\\ 0& 0& 1& 0& 0& 0 \\ 0& -1& 1& 1& 0& -1 \\ 0& -1& 1& 1& -1& 0 \\ 0& 0& 0& 1& -1& 0\end{array}\right)$,
$\left(\begin{array}{cccccc}-1& 0& 1& 0& 0& 0 \\ 0& 1& 0& 0& 0&
0\\0& 0& 1& 0& 0& 0 \\ 0& 0& 0& 1& 0& 0\\0& 0& 0& 0& 1& 0 \\ 0& 0&
0& 0& 0& 1 \end{array}\right)$,\\$\left(\begin{array}{cccccc}
-1& 1& 1& -1& 0& 0 \\ -1& 0& 2& -1& 0& 0 \\ -2& 1& 2& -1& 0& 0 \\-3&
0& 3& -1& 0& 0 \\ -2& 0& 2& 0& -1& 0 \\ -1& 0& 1& 0& 0& -1
\end{array}\right)$, $\left(\begin{array}{cccccc}-1& 1& 1& -1&
0& 0 \\ -1& 1& 1& -1& 0& 1
\\ -2& 1& 2& -1& 0& 0\\-3& 1& 2& -1& 0& 1 \\ -2& 1& 1& -1& 1& 0 \\ -1& 0& 1& -1& 1& 0\end{array}\right)$,
$\left(\begin{array}{cccccc} -1& 1& 1& -1& 0& 0 \\ -1& 1& 1& -1& 1&
-1 \\ -2& 1& 2& -1& 0& 0\\-3& 1& 2& -1& 1& -1 \\ -2& 1& 1& -1& 1& 0 \\
-1& 1& 0& 0& 0& 0\end{array}\right)$,\\
$\left(\begin{array}{cccccc}-1& 1& 1& -1& 0& 0 \\ -1& 1& 1& 0& -1& 0
\\ -2& 1& 2& -1& 0& 0\\-3& 1& 2& 0& -1& 0 \\ -2& 0& 2& 0& -1& 0 \\
       -1& 0& 1& 0& -1& 1 \end{array}\right)$,
$\left(\begin{array}{cccccc}-1& 1& 1& -1& 0& 0 \\ 0& 1& 1& -1& -1& 1
\\ 0& 2& 1& -2& 0& 0 \\ 0& 2& 1& -2& -1& 1 \\ 0& 1& 1& -1& -1& 0 \\ 0& 0& 1& -1& 0& 0\end{array}\right)$,
$\left(\begin{array}{cccccc}-1& 1& 1& -1& 0& 0 \\ 0& 1& 1& -1& 0& -1
\\ 0& 2& 1& -2& 0& 0 \\ 0& 2& 1& -2& 0& -1\\0& 1& 1& -1& -1& 0 \\ 0& 1& 0& 0& -1& 0
\end{array}\right)$,\\
$\left(\begin{array}{cccccc}1& 0& 0& 0& 0& 0 \\ 0& 0& 1& -1& 1&
0\\0& 0& 1& 0& 0& 0 \\ 0& -1& 1& 0& 1& 0\\0& 0& 0& 0& 1& 0 \\ 0& 0&
0& 0& 1& -1\end{array}\right)$, $\left(\begin{array}{cccccc}-1&
0& 1& 0& 0& 0 \\ -1& 0& 2& -1& 0& 0
\\ -2& 1& 2& -1& 0& 0 \\ -3& 0& 3& -1& 0& 0 \\-2& 0& 2& 0& -1& 0 \\ -1& 0& 1& 0& 0& -1\end{array}\right).$

The character table of $G^{s_4}$:\\
\noindent \begin{tabular}{c|cccccccccccccccccccc}
  & & & & & & & & & &10 & & & & & & & & & &20\\\hline
  $\chi _4^ {(1)}$&1&1&1& 1&1& 1&1&1&1&1&1&1&1&1&1&1&1&1&1
&1\\

  $\chi _4^ {(2)}2$&1&1&1&-1&1&-1 &-1&1 &-1&1&1 &-1&1 &-1 &-1&1 &-1 &-1&1
&-1\\

$\chi _4^ {(3)}$&1&1 &-1& 1 &-1& 1 &-1&1 &-1 &-1 &-1 &-1&1&1&1&1&1&1&1
&-1\\

$\chi _4^ {(4)}$&1&1 &-1&-1 &-1&-1&1&1&1 &-1 &-1&1&1 &-1 &-1&1 &-1 &-1&1
&1\\

$\chi _4^ {(5)}$& 1 &-1&A&-A &-A& A &-1&1&1&A &-A &-1 &-1 &-A&A&1 &-A&A &-1
&1\\

$\chi _4^ {(6)}$&1 &-1 &-A& A&A&-A &-1&1&1 &-A&A &-1 &-1&A &-A&1&A &-A &-1
&1\\

$\chi _4^ {(7)}$&1 &-1 &-A&-A&A& A&1&1 &-1 &-A&A&1 &-1 &-A&A&1 &-A&A &-1
&-1\\

$\chi _4^ {(8)}$&1 &-1&A& A &-A&-A&1&1 &-1&A &-A&1 &-1&A &-A&1&A &-A &-1
&-1\\

$\chi _4^ {(9)}$& 2 &-1&.&-1&.&-1&. &-1&.&.&.&.&2&2&2&2&2&2&2
&.\\

$\chi _4^ {({10})}$&2 &-1&.& 1&.& 1&. &-1&.&.&.&.&2 &-2 &-2&2 &-2 &-2&2
&.\\

$\chi _4^ {({11})}$&2&1&.& A&.&-A&. &-1&.&.&.&. &-2&B &-B&2&B &-B &-2
&.\\

$\chi _4^ {({12})}$&2&1&.&-A&.& A&. &-1&.&.&.&. &-2 &-B&B&2 &-B&B &-2
&.\\

$\chi _4^ {({13})}$&3&.&1& .&1& .&1&. &-1 &-1 &-1 &-1 &-1 &-1 &-1 &-1&3&3&3
&1\\

$\chi _4^ {({14})}$&3&. &-1& . &-1& . &-1&.&1&1&1&1 &-1 &-1 &-1 &-1&3&3&3
&-1\\

$\chi _4^ {({15})}$&3&. &-1& . &-1& .&1&. &-1&1&1 &-1 &-1&1&1 &-1 &-3 &-3&3
&1\\

$\chi _4^ {({16})}$& 3&.&1& .&1& . &-1&.&1 &-1 &-1&1 &-1&1&1 &-1 &-3 &-3& 3
&-1\\

$\chi _4^ {({17})}$& 3&.&A& . &-A& . &-1&. &-1 &-A&A&1&1&A &-A &-1&C &-C &-3
&1\\

$\chi _4^ {({18})}$& 3&. &-A& .&A& . &-1&. &-1&A &-A&1&1 &-A&A &-1 &-C&C &-3
&1\\

$\chi _4^ {({19})}$& 3&. &-A& .&A& .&1&.&1&A &-A &-1&1&A &-A &-1&C &-C &-3
&-1\\

$\chi _4^ {({20})}$& 3&.&A& . &-A& .&1&.&1 &-A&A &-1&1 &-A&A &-1 &-C&C &-3
  &-1\\
\end{tabular}

\noindent \noindent where A = E(4)
  = ER(-1) = i,
B = -2*E(4)
  = -2*ER(-1) = -2i,
C = -3*E(4)
  = -3*ER(-1) = -3i.

The generators of $G^{s_5}$ are:\\
$\left(\begin{array}{cccccc} 0& 0& 1& 0& -1& 0 \\ 0& 1& 1& -1& -1& 1
\\ 1& 1& 1& -1& -1& 0 \\ 1& 1& 1& -1& -2& 1\\0& 1& 1& -1& -1& 0 \\ 0& 0& 1& -1& 0&
0\end{array}\right)$,
$\left(\begin{array}{cccccc}-1& 1& 1& -1& 0& 0 \\ -1& 0& 2& -1& 0& 0
\\ -2& 1& 2& -1& 0& 0 \\ -3& 0& 3& -1& 0& 0 \\ -2& 0& 2& 0&
-1& 0 \\ -1& 0& 1& 0& 0& -1\end{array}\right)$,
$\left(\begin{array}{cccccc} -1& 1& 1& -1& 0& 0 \\ 0& 1& 1& -1& -1&
1 \\ 0& 2& 1& -2& 0& 0 \\ 0& 2& 1& -2& -1& 1 \\0& 1& 1& -1& -1& 0 \\
0& 0& 1& -1& 0& 0\end{array}\right)$.

The representatives   of conjugacy classes of   $G^{s_5}$ are:\\
$\left(\begin{array}{cccccc}1& 0& 0& 0& 0& 0 \\ 0& 1& 0& 0& 0& 0\\0&
0& 1& 0& 0& 0 \\ 0& 0& 0& 1& 0& 0 \\ 0& 0& 0& 0& 1& 0 \\ 0& 0& 0& 0&
0& 1 \end{array}\right)$,
 $\left(\begin{array}{cccccc} -1& 1& 1& -1& 0& 0 \\ -1& 1& 1& 0& -1& 0 \\ -2& 1& 2& -1& 0& 0 \\-3& 1& 2& 0& -1& 0 \\ -2& 0& 2& 0& -1& 0 \\ -1& 0& 1& 0& -1& 1 \end{array}\right)$,
 $\left(\begin{array}{cccccc}-1& 1& 1& -1& 0& 0 \\ 0& 1& 1& -1& -1& 1 \\ 0& 2& 1& -2& 0& 0 \\ 0& 2& 1& -2&
-1& 1\\0& 1& 1& -1& -1& 0 \\ 0& 0& 1& -1& 0&
0\end{array}\right)$,\\
 $\left(\begin{array}{cccccc}-1& 1& 1& -1&
0& 0 \\ 0& 1& 1& -1& 0& -1
\\ 0& 2& 1& -2& 0& 0 \\ 0& 2& 1& -2& 0& -1\\0& 1& 1& -1& -1& 0 \\ 0& 1& 0& 0& -1& 0\end{array}\right)$,
 $\left(\begin{array}{cccccc}0& 0& 1& 0& -1& 0 \\ -1& 0& 2& -1& 0& 0 \\ -1& 0& 2& 0& -1& 0 \\-2& -1& 3& 0& -1& 0 \\ -2& 0& 2& 0& -1& 0 \\ -1& 0& 1& 0& 0&
 -1\end{array}\right)$,
$\left(\begin{array}{cccccc}0& 0& 1& 0& -1& 0 \\ -1& 1& 1& 0& -1& 0
\\ -1& 0& 2& 0& -1& 0 \\ -2& 0& 2& 1& -2& 0\\-2& 0& 2& 0& -1& 0 \\ -1& 0& 1& 0& -1&
1\end{array}\right)$,\\ $\left(\begin{array}{cccccc} 0& 0& 1& 0&
-1& 0 \\ 0& 1& 1& -1& -1& 1
\\ 1& 1& 1& -1& -1& 0 \\ 1& 1& 1& -1& -2& 1\\0& 1& 1& -1& -1& 0 \\ 0& 0& 1& -1& 0& 0\end{array}\right)$,
$\left(\begin{array}{cccccc} 0& 0& 1& 0& -1& 0 \\ 0& 1& 1& -1& 0& -1
\\ 1& 1& 1& -1& -1& 0\\1& 1& 1& -1& -1& -1 \\ 0& 1& 1& -1& -1& 0 \\ 0& 1& 0& 0& -1& 0 \end{array}\right)$,
$\left(\begin{array}{cccccc}0& 1& 0& -1& 1& 0 \\ 0& 0& 1& -1& 1& 0
\\ -1& 1& 1& -1& 1& 0 \\ -1& 0& 1& -1& 2& 0\\0& 0& 0& 0& 1& 0 \\ 0& 0& 0& 0& 1& -1
\end{array}\right)$,\\
$\left(\begin{array}{cccccc}0& 1& 0& -1& 1& 0 \\ 0& 1& 0& 0& 0& 0
\\-1& 1& 1& -1& 1& 0 \\ -1& 1& 0& 0& 1& 0\\0& 0& 0& 0& 1& 0 \\ 0& 0& 0& 0& 0& 1\end{array}\right)$,
$\left(\begin{array}{cccccc}0& 1& 0& -1& 1& 0 \\ 1& 1& 0& -1& 0& 1
\\ 1& 2& 0& -2& 1& 0 \\ 2& 2& -1& -2& 1& 1 \\ 2& 1& -1& -1& 1& 0 \\ 1& 0& 0& -1& 1& 0\end{array}\right)$,
$\left(\begin{array}{cccccc} 0& 1& 0& -1& 1& 0 \\ 1& 1& 0& -1& 1& -1
\\ 1& 2& 0& -2& 1& 0 \\ 2& 2& -1& -2& 2& -1 \\2& 1& -1& -1& 1& 0 \\ 1& 1& -1& 0& 0& 0
\end{array}\right)$,\\
$\left(\begin{array}{cccccc} 1& 0& 0& 0& 0& 0 \\ 0& 0& 1& -1& 1& 0
\\ 0& 0& 1& 0& 0& 0 \\ 0& -1& 1& 0& 1& 0\\0& 0& 0& 0& 1& 0 \\ 0& 0& 0& 0& 1&
-1\end{array}\right)$, $\left(\begin{array}{cccccc} -1& 1& 1&
-1& 0& 0 \\ -1& 0& 2& -1& 0&
0 \\ -2& 1& 2& -1& 0& 0\\-3& 0& 3& -1& 0& 0 \\ -2& 0& 2& 0& -1& 0 \\
-1& 0& 1& 0& 0& -1\end{array}\right)$,
$\left(\begin{array}{cccccc} 1& 0& 0& 0& 0& 0 \\ 1& 1& 0& -1& 0& 1
\\ 2& 1& 0& -1& 0& 0 \\ 3& 1& -1& -1& 0& 1\\2& 1& -1& -1& 1& 0 \\ 1& 0& 0& -1& 1& 0
\end{array}\right)$,\\
$\left(\begin{array}{cccccc}1& 0& 0& 0& 0& 0 \\ 1& 1& 0& -1& 1& -1
\\ 2& 1& 0& -1& 0& 0 \\ 3& 1& -1& -1& 1& -1\\2& 1& -1& -1& 1& 0 \\ 1& 1& -1& 0& 0& 0 \end{array}\right).$

The character table of $G^{s_5}$:\\
\noindent \begin{tabular}{c|cccccccccccccccc}
  & & & & & & & & & & 10& & & & & &\\\hline
 $\chi _5^ {(1)}$     &1 &1 &1 &1 &1 &1 &1 &1 &1 &1 &1 &1 &1 &1 &1 &1\\

  $\chi_5^ {(2)}$     &1&-1&-1&-1 &1 &1 &1 &1&-1&-1&-1&-1 &1&-1 &1 &1\\

  $\chi_ 5^ {(3)}$     &1&-1&-1&-1&-1&-1&-1&-1 &1 &1 &1 &1 &1&-1 &1 &1\\

  $\chi_ 5^ {(4)}$     &1&-1 &1 &1 &1 &1&-1&-1&-1&-1 &1 &1 &1&-1&-1&-1\\

  $\chi _5^ {(5)}$     &1&-1 &1 &1&-1&-1 &1 &1 &1 &1&-1&-1 &1&-1&-1&-1\\

  $\chi _5^ {(6)}$     &1 &1&-1&-1&-1&-1 &1 &1&-1&-1 &1 &1 &1 &1&-1&-1\\

  $\chi _5^ {(7)}$     &1 &1&-1&-1 &1 &1&-1&-1 &1 &1&-1&-1 &1 &1&-1&-1\\

  $\chi _5^ {(8)}$     &1 &1 &1 &1&-1&-1&-1&-1&-1&-1&-1&-1 &1 &1 &1 &1\\

  $\chi _5^ {(9)}$     &1 &1 &A&-A &1&-1&-A &A &1&-1&-A &A&-1&-1 &A&-A\\

  $\chi _5^ {({10})}$    &1 &1&-A &A &1&-1 &A&-A &1&-1 &A&-A&-1&-1&-A &A\\

  $\chi_ 5^ {({11})}$    &1 &1 &A&-A&-1 &1 &A&-A&-1 &1 &A&-A&-1&-1 &A&-A\\

  $\chi _5^ {({12})}$    &1 &1&-A &A&-1 &1&-A &A&-1 &1&-A &A&-1&-1&-A &A\\

  $\chi_ 5^ {({13})}$    &1&-1 &A&-A&-1 &1&-A &A &1&-1 &A&-A&-1 &1&-A &A\\

  $\chi _5^ {({14})}$    &1&-1&-A &A&-1 &1 &A&-A &1&-1&-A &A&-1 &1 &A&-A\\

  $\chi _5^ {({15})}$    &1&-1 &A&-A &1&-1 &A&-A&-1 &1&-A &A&-1 &1&-A &A\\

  $\chi _5^ {({16})}$    &1&-1&-A &A &1&-1&-A &A&-1 &1 &A&-A&-1 &1 &A&-A\\

  \end{tabular}

\noindent \noindent where A = -E(4)
  = -ER(-1) = -i.

The generators of $G^{s_6}$ are:\\
$\left(\begin{array}{cccccc}-1& 1& 1& -1& 0& 0
\\ 0& 1& 1& -2& 1& 0
\\ 0& 2& 1& -2& 0& 0 \\ 0& 2& 1& -3& 1& 0\\ 0& 2& 0& -2& 1& 0 \\ 0& 1& 0& -1& 1& -1 \end{array}\right)$,
 $\left(\begin{array}{cccccc} -1& 0& 1& 0& 0& 0 \\ 0& 1& 0& 0& 0& 0 \\ 0& 0& 1& 0& 0& 0 \\ 0& 0& 0& 1& 0& 0 \\0& 0& 0& 0& 1& 0 \\ 0& 0& 0& 0& 0& 1 \end{array}\right)$,
$\left(\begin{array}{cccccc} 1& 1& 0& -1& 0& 0 \\ 0& 2& 0& -1& 0& 0 \\ 0& 2& 1& -2& 0& 0 \\
0& 3& 0& -2& 0& 0\\0& 2& 0& -2& 1& 0 \\ 0& 1& 0& -1& 0& 1
\end{array}\right)$,\\
 $\left(\begin{array}{cccccc} 1& 0& 0& 0& 0& 0 \\ 1& 0& 1& -1& 0& 0 \\ 2& 1& 0& -1& 0& 0 \\ 3& 0& 0& -1& 0& 0 \\2& 0& 0& 0& -1& 0 \\ 1& 0& 0& 0& 0& -1 \end{array}\right)$,
$\left(\begin{array}{cccccc} 1& 0& 0& 0& 0& 0 \\ 0& 1& 0& 0& 0& 0 \\ 1& 1& 0& 0& 0& -1 \\
1& 1& -1& 1& 0& -1\\1& 1& -1& 0& 1& -1 \\ 1& 1& -1& 0& 0& 0
\end{array}\right)$,
$\left(\begin{array}{cccccc} 1& 0& 0& 0& 0& 0 \\ 0& 1& 0& 0& 0& 0 \\
0& 0& 1& 0& 0& 0 \\ 0& 0& 0& 1& 0& 0\\0& 0& 0& 0& 1& 0 \\ 0& 0& 0&
0& 1& -1 \end{array}\right)$,\\
$\left(\begin{array}{cccccc} 1& 0& -1& 0& 0& 1 \\ 0& 0& -1& 0& 0& 1 \\ 0& 0& -1& 0& 0& 2 \\
0& 0& -1& 0& -1& 3\\0& -1& -1& 1& -1& 2 \\ 0& -1& 0& 0& 0& 1
\end{array}\right).$

The representatives   of conjugacy classes of   $G^{s_6}$ are:\\
$\left(\begin{array}{cccccc} 1& 0& 0& 0& 0& 0 \\ 0& 1& 0& 0& 0& 0 \\ 0& 0& 1& 0& 0& 0 \\
0& 0& 0& 1& 0& 0\\0& 0& 0& 0& 1& 0 \\ 0& 0& 0& 0& 0& 1
\end{array}\right)$,
  $\left(\begin{array}{cccccc} -1& 1& 1& -1& 0& 0 \\ 0& 1& 1& -2& 1& 0 \\ 0& 2& 1& -2& 0& 0 \\ 0& 2& 1& -3& 1&
0 \\ 0& 2& 0& -2& 1& 0 \\ 0& 1& 0& -1& 1& -1
\end{array}\right)$,
  $\left(\begin{array}{cccccc} -1& 1& 1& -1& 0& 0 \\ -1& 1& 1& -1& 0& 1 \\ -2& 1& 2& -1& 0& 0 \\-3& 1& 2& -1& 0& 1 \\ -2& 1& 1& -1& 1& 0 \\ -1& 1& 0& 0& 0& 0
  \end{array}\right)$,\\
  $\left(\begin{array}{cccccc} -1& 0& 1& 0& 0& 0 \\ -1& 1& 1& -1& 0& 1 \\ -2& 1& 2& -1& 0& 0 \\ -3& 1& 2& -1&
0& 1 \\ -2& 1& 1& -1& 1& 0 \\ -1& 0& 1& -1& 1& 0
\end{array}\right)$,
  $\left(\begin{array}{cccccc} -1& 0& 1& 0& 0& 0 \\ 0& 1& 0& 0& 0& 0 \\ -1& 1& 1& 0& 0& -1 \\ -1& 1& 0& 1& 0&
-1\\ -1& 1& 0& 0& 1& -1 \\ -1& 1& 0& 0& 0& 0
\end{array}\right)$,
  $\left(\begin{array}{cccccc} 1& 1& 0& -1& 0& 0 \\ 0& 1& 1& -2& 1& 0 \\ 1& 1& 1& -2& 0& 1 \\ 1& 1& 1& -3& 1&
1\\ 1& 1& 0& -2& 1& 1 \\ 1& 0& 0& -1& 1& 0
\end{array}\right)$,\\
  $\left(\begin{array}{cccccc} -1& 0& 1& 0& 0& 0 \\ 0& 1& 0& 0& 0& 0 \\ 0& 0& 1& 0& 0& 0 \\ 0& 0& 0& 1& 0& 0
\\ 0& 0& 0& 0& 1& 0 \\ 0& 0& 0& 0& 0& 1 \end{array}\right)$,
  $\left(\begin{array}{cccccc} 1& 1& 0& -1& 0& 0 \\ 0& 1& 1& -2& 1& 0 \\ 0& 2& 1& -2& 0& 0 \\ 0& 2& 1& -3& 1&
0 \\0& 2& 0& -2& 1& 0 \\ 0& 1& 0& -1& 1& -1 \end{array}\right)$,
  $\left(\begin{array}{cccccc} -1& 0& 1& 0& 0& 0 \\ -1& 1& 1& -1& 0& 1 \\ -2& 1& 2& -1& 0& 0 \\ -3& 1& 2& -1&
0& 1 \\-2& 1& 1& -1& 1& 0 \\ -1& 1& 0& 0& 0& 0
\end{array}\right)$,\\
  $\left(\begin{array}{cccccc} 0& 0& 0& 0& -1& 1 \\ 0& 0& 0& 1& -2& 1 \\ 0& -1& 0& 1& -2& 2 \\ 0& 0& 0& 1& -3&
3\\ 0& 0& 1& 0& -2& 2 \\ -1& 0& 1& 0& -1& 1 \end{array}\right)$,
  $\left(\begin{array}{cccccc} 0& -1& 0& 1& 0& -1 \\ 0& -1& 1& 0& 0& -1 \\ 0& -1& 0& 1& 0& -2\\0& -1& 1& 0& 1& -3 \\ 0& 0& 1& 0& 0& -2 \\ 1& 0& 0& 0& 0& -1
  \end{array}\right)$,
  $\left(\begin{array}{cccccc} 0& 0& 0& 0& -1& 1 \\ 0& 0& 0& 1& -2& 1 \\ -1& 0& 0& 1& -2& 1 \\ -1& 1& 0& 1& -
3& 2 \\ -1& 1& 1& 0& -2& 1 \\ 0& 1& 0& 0& -1& 0
\end{array}\right)$,\\
  $\left(\begin{array}{cccccc} 1& 0& -1& 1& -1& 0 \\ 0& 0& 0& 1& -2& 1 \\ 1& -1& -1& 2& -2& 1\\1& 0& -1& 2& -3& 2 \\ 1& 0& 0& 1& -2& 1 \\ 1& 0& 0& 0& -1& 1 \end{array}\right)$,
  $\left(\begin{array}{cccccc} -1& -1& 1& 0& 0& 0 \\ 0& -1& 1& 0& 0& -1 \\ -1& -1& 1& 0& 0& -1\\-1& -1& 2& -1& 1& -2 \\ -1& 0& 2& -1& 0& -1 \\ -1& 0& 1& 0& 0& -1 \end{array}\right)$,
  $\left(\begin{array}{cccccc} 0& 0& 0& 0& -1& 1 \\ 0& 0& 0& 1& -2& 1 \\ 0& -1& 0& 1& -2& 2 \\ 0& 0& 0& 1& -3&
3 \\0& 0& 1& 0& -2& 2 \\ 1& 0& 0& 0& -1& 1
\end{array}\right)$,\\
  $\left(\begin{array}{cccccc} 0& -1& 0& 1& 0& -1 \\ 0& -1& 1& 0& 0& -1 \\ 0& -1& 0& 1& 0& -2 \\0& -1& 1& 0& 1& -3 \\ 0& 0& 1& 0& 0& -2 \\ -1& 0& 1& 0& 0& -1
  \end{array}\right)$,
  $\left(\begin{array}{cccccc} 1& 0& 0& 0& 0& 0 \\ 1& 0& 1& -1& 0& 0 \\ 2& 1& 0& -1& 0& 0 \\ 3& 0& 0& -1& 0& 0
\\ 2& 0& 0& 0& -1& 0 \\ 1& 0& 0& 0& 0& -1 \end{array}\right)$,
  $\left(\begin{array}{cccccc} -1& 1& 1& -1& 0& 0 \\ -1& 1& 1& 0& -1& 0 \\ -2& 1& 2& -1& 0& 0\\-3& 1& 2& 0& -1& 0 \\ -2& 0& 2& 0& -1& 0 \\ -1& 0& 1& 0& -1& 1
  \end{array}\right)$,\\
  $\left(\begin{array}{cccccc} -1& 0& 1& 0& 0& 0 \\ -1& 0& 2& -1& 0& 0 \\ -2& 1& 2& -1& 0& 0 \\ -3& 0& 3& -1&
0& 0 \\-2& 0& 2& 0& -1& 0 \\ -1& 0& 1& 0& -1& 1
\end{array}\right)$,
  $\left(\begin{array}{cccccc} -1& 0& 1& 0& 0& 0 \\ -1& 0& 2& -1& 0& 0 \\ -1& 0& 2& -1& 0& 1\\-2& -1& 3& -1& 0& 1 \\ -1& -1& 2& 0& -1& 1 \\ 0& -1& 1& 0& 0& 0
  \end{array}\right)$,
  $\left(\begin{array}{cccccc} 1& 1& 0& -1& 0& 0 \\ 1& 1& 0& 0& -1& 0 \\ 1& 2& 0& -1& 0& -1 \\ 2& 2& -1& 0& -
1& -1\\1& 1& 0& 0& -1& -1 \\ 0& 1& 0& 0& -1& 0
\end{array}\right)$,\\
  $\left(\begin{array}{cccccc} -1& 0& 1& 0& 0& 0 \\ -1& 0& 2& -1& 0& 0 \\ -2& 1& 2& -1& 0& 0 \\ -3& 0& 3& -1&
0& 0\\-2& 0& 2& 0& -1& 0 \\
       -1& 0& 1& 0& 0& -1 \end{array}\right)$,
  $\left(\begin{array}{cccccc} 1& 1& 0& -1& 0& 0 \\ 1& 1& 0& 0& -1& 0 \\ 2& 1& 0& -1& 0& 0 \\ 3& 1& -1& 0& -1&
0 \\2& 0& 0& 0& -1& 0 \\ 1& 0& 0& 0& -1& 1 \end{array}\right)$,
  $\left(\begin{array}{cccccc} -1& 0& 1& 0& 0& 0 \\ 0& 0& 1& 0& 0& -1 \\ 0& 0& 1& 0& 0& 0 \\ 0& -1& 1& 1& 0&
-1\\0& -1& 1& 1& -1& 0 \\ 0& -1& 1& 0& 0& 0
\end{array}\right)$,\\
  $\left(\begin{array}{cccccc} 0& 1& 0& 0& -1& 1 \\ 0& 1& 1& -1& -1& 1 \\ 0& 1& 1& -1& -1& 2 \\ 0& 1& 1& -1&
-2& 3\\0& 0& 1& 0& -2& 2 \\ 1& 0& 0& 0& -1& 1
\end{array}\right).$

 The character table of $G^{s_6}$:\\
\noindent
\begin{tabular}{c|ccccccccccccccccccccccccc}
  & & & & & & & & & &10 & & & & & & & & & &20 & & & & &\\\hline
$\chi_6^ {(1)}$&1& 1&1&1&1&1&1&1&1&1&1& 1&1&1&1&1&1&1&1&1&1&1&1&1&1\\

$\chi_6^ {(2)}$& 1& 1&1&1&1&1&-1&-1&-1&1&1& 1&1&1&-1&-1&-1&-1&-1&-1
&-1&1&1&1
&1\\

$\chi_6^ {(3)}$& 1& 1&1&1&1&1&-1&-1&-1&1&1& 1&1&1&-1&-1&1&1&1&1&1&-1
&-1&-1
&-1\\

$\chi_6^ {(4)}$& 1& 1&1&1&1&1&1&1&1&1&1&
1&1&1&1&1&-1&-1&-1&-1&-1&-1&-1 &-1&
-1\\

$\chi_6^ {(5)}$& 2& 2&2&2&2&2&-2&-2&-2&-1&-1&-1&-1&-1&1&1&.&.&.&.&.&.&.&. &.\\

$\chi_6^ {()}6$& 2& 2&2&2&2&2&2&2&2&-1&-1&-1&-1&-1&-1&-1&.&.&.&.&.&. &.&.&.\\

$\chi_6^ {(7)}$&2& 2&2&2&-1&-1&.&.&.&2&2&-1&-1&-1&.&.&-2&-2&-2&1&1&. &.&.&.\\
\end{tabular}

The generators of $G^{s_7}$ are:\\
 $\left(\begin{array}{cccccc} 1& 0& -1& 0& 1& -1
\\ 1& 1& -1& -1&1& 0\\ 1& 0& -1&-1& 2& -1 \\ 2& 1& -1&-2& 2&  -1 \\ 1& 1& 0& -2& 2&
-1 \\ 0& 1& 0& -1& 1& -1
\end{array}\right)$, $\left(\begin{array}{cccccc} -1& 0& 0& 1& -1& 0
\\ -1& -1& 1& 1& -1&
0 \\ -2& 0& 1& 1& -1& 0\\-2& 0& 2& 1& -2& 0 \\ -2& 0& 2& 0& -1& 0 \\
-1& 0& 1& 0& 0& -1 \end{array}\right)$,
 $\left(\begin{array}{cccccc}
1& 0& -1& 0& 1& -1 \\ 0& 1& -1& 0& 1& -1 \\ 0& 0& -1& 0& 2& -2\\0&
0& -2& 1& 2& -2 \\ 0& 0& -1& 0& 2& -1 \\ 0& 0& -1& 0& 1& 0
\end{array}\right)$,\\
$\left(\begin{array}{cccccc} 1& 0& 0& 0& 0& 0
\\ 1& 0& -1& 1& -1& 1
\\ 1& -1& 0& 1& -1& 1 \\ 1& -1& -1& 2& -1& 1 \\1& -1& -1& 2& -1& 0 \\ 0& 0& 0& 1& -1& 0 \end{array}\right).$

The representatives   of conjugacy classes of   $G^{s_7}$ are:\\
$\left(\begin{array}{cccccc} 1& 0& 0& 0& 0& 0 \\ 0& 1& 0& 0& 0& 0 \\ 0& 0& 1& 0& 0& 0 \\
0& 0& 0& 1& 0& 0\\0& 0& 0& 0& 1& 0 \\ 0& 0& 0& 0& 0& 1
\end{array}\right)$,
  $\left(\begin{array}{cccccc} -1& 0& 0& 0& 0& 1 \\ -1& -1& 0& 0& 1& 0 \\ -1& -1& -1& 1& 0& 1\\-1& -1& -1& 0& 1& 1 \\ -1& -1& 0& 0& 0& 1 \\ 0& -1& 0& 0& 0& 1
  \end{array}\right)$,
  $\left(\begin{array}{cccccc} -1& 0& 0& 0& 0& 1 \\ -1& 0& 0& 0& 1& 0 \\ -1& 0& -1& 1& 0& 1 \\ -2& 0& -1& 1&
1& 0\\-1& 1& -1& 0& 1& 0 \\ 0& 0& -1& 0& 1& 0
\end{array}\right)$,\\
  $\left(\begin{array}{cccccc} -1& 0& 0& 0& 0& 1 \\ -1& 0& 0& 0& 1& 0 \\ -1& 0& -1& 1& 0& 1 \\ -2& 0& 0& 1& 0&
1\\ -1& 1& 0& 0& 0& 1 \\ 0& 0& 0& 0& 0& 1 \end{array}\right)$,
  $\left(\begin{array}{cccccc} -1& 0& 0& 0& 0& 1 \\ 0& 0& 0& -1& 1& 1 \\ 0& 0& -1& 0& 0& 2 \\ 0& 0& -1& -1& 1&
2\\ 0& 1& -1& -1& 1& 1 \\ 0& 0& -1& 0& 1& 0 \end{array}\right)$,
  $\left(\begin{array}{cccccc} -1& 0& 0& 0& 0& 1 \\ 0& 0& 0& -1& 1& 1 \\ 0& 0& -1& 0& 0& 2 \\ 0& 0& 0& -1& 0&
3\\ 0& 1& 0& -1& 0& 2 \\ 0& 0& 0& 0& 0& 1 \end{array}\right)$,\\
  $\left(\begin{array}{cccccc} 0& -1& 0& 0& 1& 0 \\ -1& -1& 0& 0& 1& 0 \\ 0& -2& -1& 1& 1& 0 \\-1& -2& -1& 1& 2& -1 \\ -1& -2& 0& 1& 1& -1 \\ -1& -1& 0& 1& 0& 0 \end{array}\right)$,
  $\left(\begin{array}{cccccc} 0& -1& 0& 0& 1& 0 \\ -1& -1& 0& 0& 1& 0 \\ 0& -2& -1& 1& 1& 0 \\ 0& -2& -1& 0&
2& 0\\0& -2& 0& 0& 1& 0 \\ 0& -1& 0& 0& 0& 1
\end{array}\right)$,
  $\left(\begin{array}{cccccc} 0& -1& 0& 0& 1& 0 \\ -1& 0& 0& 0& 1& 0 \\ 0& -1& -1& 1& 1& 0\\-1& -1& -1& 1& 2& -1 \\ 0& 0& -1& 0& 2& -1 \\ 0& 0& -1& 0& 1& 0
  \end{array}\right)$,\\
  $\left(\begin{array}{cccccc} 0& -1& 0& 0& 1& 0 \\ -1& 0& 0& 0& 1& 0 \\ 0& -1& -1& 1& 1& 0 \\ -1& -1& 0& 1&
1& 0\\ 0& 0& 0& 0& 1& 0 \\ 0& 0& 0& 0& 0& 1 \end{array}\right)$,
  $\left(\begin{array}{cccccc} 0& -1& 0& 0& 1& 0 \\ 0& 0& 0& -1& 1& 1 \\ 1& -1& -1& 0& 1& 1 \\ 1& -1& -1& -1&
2& 1\\ 1& 0& -1& -1& 2& 0 \\ 0& 0& -1& 0& 1& 0
\end{array}\right)$,
  $\left(\begin{array}{cccccc} 0& -1& 0& 0& 1& 0 \\ 0& 0& 0& -1& 1& 1 \\ 1& -1& -1& 0& 1& 1 \\ 1& -1& 0& -1&
1& 2 \\ 1& 0& 0& -1& 1& 1 \\ 0& 0& 0& 0& 0& 1
\end{array}\right)$,\\
  $\left(\begin{array}{cccccc} 0& 1& -1& 0& 0& 0 \\ 0& 1& -1& 0& 1& -1 \\ -1& 1& -1& 0& 1& -1\\-1& 1& -2& 1& 1& -1 \\ -1& 1& -1& 0& 1& 0 \\ 0& 0& -1& 0& 1& 0 \end{array}\right)$,
  $\left(\begin{array}{cccccc} 0& 1& -1& 0& 0& 0 \\ 0& 1& -1& 0& 1& -1 \\ -1& 1& -1& 0& 1& -1 \\-1& 2& -1& 0& 1& -2 \\ -1& 2& 0& -1& 1& -1 \\ 0& 1& 0& -1& 1& -1 \end{array}\right)$,
  $\left(\begin{array}{cccccc} 0& 1& -1& 0& 0& 0 \\ 1& 1& -1& -1& 1& 0 \\ 0& 1& -1& -1& 1& 0 \\ 1& 2& -1& -2&
1& 0\\ 0& 2& 0& -2& 1& 0 \\ 0& 1& 0& -1& 1& -1
\end{array}\right)$,\\
  $\left(\begin{array}{cccccc} 0& 1& 0& 0& -1& 1 \\ 0& 1& 0& 0& 0& 0 \\ -1& 1& 1& 0& -1& 1 \\ -1& 1& 0& 1& -1&
1\\ -1& 1& 0& 0& 0& 1 \\ 0& 0& 0& 0& 0& 1 \end{array}\right)$,
  $\left(\begin{array}{cccccc} 1& 0& -1& 0& 1& -1 \\ 0& 1& -1& 0& 1& -1 \\ 0& 0& -1& 0& 2& -2\\0& 0& -2& 1& 2& -2 \\ 0& 0& -1& 0& 2& -1 \\ 0& 0& -1& 0& 1& 0
  \end{array}\right)$,
  $\left(\begin{array}{cccccc} 1& 0& -1& 0& 1& -1 \\ 0& 1& -1& 0& 1& -1 \\ 0& 0& -1& 0& 2& -2\\ 0& 1& -1& 0& 2& -3 \\ 0& 1& 0& -1& 2& -2 \\ 0& 1& 0& -1& 1& -1
  \end{array}\right)$,\\
  $\left(\begin{array}{cccccc} 1& 0& -1& 0& 1& -1 \\ 1& 1& -1& -1& 1& 0 \\ 1& 0& -1& -1& 2& -1\\2& 1& -1& -2& 2& -1 \\ 1& 1& 0& -2& 2& -1 \\ 0& 1& 0& -1& 1& -1
  \end{array}\right)$,
  $\left(\begin{array}{cccccc} -1& 0& 0& 0& 0& 1 \\ -1& -1& 0& 0& 1& 0 \\ -1&
-1& -1& 1& 0& 1\\-2& -1& -1& 1& 1& 0 \\ -2& -1& 0& 1& 0& 0 \\ -1& -1& 0& 1& 0& 0 \end{array}\right).$

The character table of $G^{s_7}$:\\
\noindent \begin{tabular}{c|cccccccccccccccccccc}
  & & & & & & & & & &10 & & & & & & & & & &20\\\hline
$\chi_7^ {(1)}$&1&1&1&1&1&1&1&1&1&1&1&1&1&1&1&1&1&1&1
&1\\

$\chi_7^ {(2)}$& 1&-1&-1&1&1&-1&-1&1&1&-1&-1&1&1&-1&1&-1&-1 &1&-1
&1\\

$\chi_7^ {(3)}$&1&1&-1&1&1&-1&-1&1&-1&1&1&-1&-1&1&-1&1&-1&1 &-1
&-1\\

$\chi_7^ {(4)}$& 1&1&-1&-1&-1&-1& 1&1&-1&-1&-1&-1&1&1&1&1&1 & 1&1
&1\\

$\chi_7^ {(5)}$&1&-1&-1&-1&-1&-1&1&1&1&1&1&1&-1&-1&-1&-1&1&1 &1
&-1\\

$\chi_7^ {(6)}$&1&1&1&-1&-1&1&-1&1&1&-1&-1&1&-1&1&-1&1&-1&1 &-1
&-1\\

$\chi_7^ {(7)}$&1&-1&1&-1&-1&1&-1&1&-1&1&1&-1&1&-1&1&-1&-1&1 &-1
&1\\

$\chi_7^ {(8)}$&1&-1&1&1&1&1&1&1&-1&-1&-1&-1&-1&-1&-1&-1&1&1 &1
&-1\\

$\chi_7^ {(9)}$&2&1&.&.&.&.&-1&-1&.&.&.&.&-2&-2&-2&-2&2&2&2
&1\\

$\chi_7^ {({10})}$&2&-1&.&.&.&.&1&-1&.&.&.&.&-2&2&-2&2&-2&2 &-2
&1\\

$\chi_7^ {({11})}$&2&-1&.&.&.&.&-1&-1&.&.&.&.&2&2&2&2&2&2&2
&-1\\

$\chi_7^ {({12})}$&2&1&.&.&.&.&1&-1&.&.&.&.&2&-2&2&-2&-2&2 &-2
&-1\\

$\chi_7^ {({13})}$&3&.&1&-1&1&-1&.&.&-1&1&-1&1&1&1&-3&-3&-1 &-1&3
&.\\

$\chi_7^ {({14})}$&3&.&-1&1&-1&1&.&.&1&-1&1&-1&1&1&-3&-3&-1 &-1&3
&.\\

$\chi_7^ {({15})}$&3&.&1&1&-1&-1&.&.&1&1&-1&-1&1&-1&-3&3&1 & -1&-3
&.\\

$\chi_7^ {({16})}$&3&.&-1&-1&1&1&.&.&-1&-1&1&1&1&-1&-3&3&1 &-1&-3
&.\\

$\chi_7^ {({17})}$&3&.&1&-1&1&-1&.&.&1&-1&1&-1&-1&-1&3&3&-1 & -1& 3
& .\\

$\chi_7^ {({18})}$&3&.&-1&1&-1&1&.&.&-1&1&-1&1&-1&-1&3&3&-1 &-1 &3
&.\\

$\chi_7^ {({19})}$&3&.&1&1&-1&-1&.&.&-1&-1&1&1&-1&1&3&-3&1 &-1& -3
& .\\

$\chi_7^ {({20})}$&3&.&-1&-1&1&1&.&.&1&1&-1&-1&-1&1&3&-3&1 &-1&-3
& .\\
\end{tabular}


The generators of $G^{s_8}$ are:\\
$\left(\begin{array}{cccccc} 1& -1& 0& 0& 0& 0
\\ 1& 0& 0& 0& -1& 0\\2& 0& 0& -1& 0& 0 \\ 2& 0& 1& -1& -1& 0\\2& 0& 0& 0& -1& 0 \\ 1& 0& 0& 0& -1& 1 \end{array}\right)$,
$\left(\begin{array}{cccccc} 1& 0& 0& 0& 0& 0 \\ 0& 1& 0& 0& 0& 0 \\
0& 0& 1& 0& 0& 0 \\ 0& 0& 0& 1& 0& 0\\0& 0& 0& 0& 1& 0 \\ 0& 0& 0&
0& 1& -1 \end{array}\right)$,
 $\left(\begin{array}{cccccc} -1& 1& 1& -1& 0& 0 \\ -1& 0& 2& -1& 0& 0 \\ -2& 1& 2& -1& 0& 0\\-3& 0& 3& -1& 0& 0 \\ -2& 0& 2& 0& -1& 0 \\ -1& 0& 1& 0& 0& -1
 \end{array}\right)$,\\
 $\left(\begin{array}{cccccc} 0& 1& 0& 0& 0& -1 \\ 0& 1& 0& 0& 0& 0 \\ -1& 1& 1& 0& 0& -1 \\ -1& 1& 0& 1& 0& -1\\-1& 1& 0& 0& 1& -1 \\ -1& 1& 0& 0& 0& 0
 \end{array}\right)$,
 $\left(\begin{array}{cccccc} -1& 1& 1& -1& 0& 0 \\ -1& 0& 2& -1& 0& 0
\\ -2& 0& 2& -1& 0& 0\\-2& 0& 3& -1& -1& 0 \\ -2& 0& 2& 0& -1& 0 \\ -1& 0& 1& 0& 0& -1 \end{array}\right).$

The representatives   of conjugacy classes of   $G^{s_8}$ are:\\
$\left(\begin{array}{cccccc} 1& 0& 0& 0& 0& 0 \\ 0& 1& 0& 0& 0& 0 \\ 0& 0& 1& 0& 0& 0 \\
0& 0& 0& 1& 0& 0  \\0& 0& 0& 0& 1& 0 \\ 0& 0& 0& 0& 0&
1\end{array}\right)$ ,
  $\left(\begin{array}{cccccc} -1& 0& 0& 0& 0& 1 \\ -1& -1& 0& 0& 1& 0 \\ -2& -1& 1& 0& 0& 1 \\ -2& -2& 0& 1&
0& 1 \\ -1& -1& 0& 0& 0& 1 \\ 0& -1& 0& 0& 0& 1
\end{array}\right)$,
  $\left(\begin{array}{cccccc} -1& 0& 0& 0& 0& 1 \\ -1& -1& 0& 0& 1& 0 \\ -1& -1& -1& 1& 0& 1\\
  -1& -1& -1& 0& 1& 1 \\ -1& -1& 0& 0& 0& 1 \\ -1& 0& 0& 0& 0& 0
  \end{array}\right)$,\\
  $\left(\begin{array}{cccccc} -1&
   0& 0& 0& 0& 1 \\ -1& -1& 0& 0& 1& 0 \\ -1& -1& -1& 1& 0& 1 \\-1& -1& -1& 0& 1& 1
   \\ -1& -1& 0& 0& 0& 1 \\ 0& -1& 0& 0& 0& 1 \end{array}\right)$,
  $\left(\begin{array}{cccccc} -1& 0& 0& 0& 0& 1 \\ -1& -1& 0& 0& 1& 0 \\ -1& 0& 0& -1& 1& 1\\-2& -1& 1& -1& 1& 1
   \\ -1& -1& 0& 0& 0& 1 \\ -1& 0& 0& 0& 0& 0 \end{array}\right)$,
  $\left(\begin{array}{cccccc} -1& 0& 0& 0& 0& 1 \\ -1& -1& 0& 0& 1& 0 \\ -1& 0& 0& -1& 1& 1 \\-2& -1& 1& -1& 1&
  1 \\ -1& -1& 0& 0& 0& 1 \\ 0& -1& 0& 0& 0& 1
  \end{array}\right)$,\\
  $\left(\begin{array}{cccccc} -1& 0& 0& 1& -1& 0 \\ 0& -1& -1& 2& -1& 0 \\ -1& -1& -1& 2& -1& 0\\ 0& -1& -2& 3&
  -2& 0 \\ 0& 0& -2& 2& -1& 0 \\ 0& 0& -1& 1& -1& 1 \end{array}\right)$,
  $\left(\begin{array}{cccccc} -1& 0& 0& 1& -1& 0 \\ 0& -1& -1& 2& -1& 0 \\ -1& -1& -1& 2& -1& 0 \\ 0& -1& -2& 3
  & -2& 0 \\ 0& 0& -2& 2& -1& 0 \\ 0& 0& -1& 1& 0& -1 \end{array}\right)$,
  $\left(\begin{array}{cccccc} -1& 0& 0& 1& -1& 0 \\ 0& -1& -1& 2& -1& 0 \\ -1& 0& -1& 2& -1& 0\\-1& -1& -2& 3&
  -1& 0 \\ 0& 0& -2& 2& -1& 0 \\ 0& 0& -1& 1& -1& 1
  \end{array}\right)$,\\
  $\left(\begin{array}{cccccc} -1& 0& 0& 1& -1& 0 \\ 0& -1& -1& 2& -1& 0 \\ -1& 0& -1& 2& -1& 0 \\-1& -1& -2& 3&
  -1& 0 \\ 0& 0& -2& 2& -1& 0 \\ 0& 0& -1& 1& 0& -1
  \end{array}\right)$,
  $\left(\begin{array}{cccccc} -1& 0& 0& 1& -1& 0 \\ 0& -1& -1& 2& -1& 0 \\ 0& 0& -1& 2& -2& 0 \\0& 0& -2& 3& -2
   0 \\ 0& 0& -2& 2& -1& 0 \\ 0& 0& -1& 1& -1& 1 \end{array}\right)$,
  $\left(\begin{array}{cccccc} -1& 0& 0& 1& -1& 0 \\ 0& -1& -1& 2& -1& 0 \\ 0& 0& -1& 2& -2& 0\\0& 0& -2& 3& -2
  & 0 \\ 0& 0& -2& 2& -1& 0 \\ 0& 0& -1& 1& 0& -1
  \end{array}\right)$,\\
  $\left(\begin{array}{cccccc} 0& -1& 0& 0& 1& 0 \\ -1& -1& 0& 0& 1& 0 \\ -1& -2& 1& 0& 1& 0 \\ -1& -3& 0& 1&
1& 0 \\ 0& -2& 0& 0& 1& 0 \\ 0& -1& 0& 0& 0& 1
\end{array}\right)$,
  $\left(\begin{array}{cccccc} 0& -1& 0& 0& 1& 0 \\ -1& -1& 0& 0& 1& 0 \\ 0& -2& -1& 1& 1& 0 \\ 0& -2& -1& 0&
2& 0 \\ 0& -2& 0& 0& 1& 0 \\0& -1& 0& 0& 0& 1
\end{array}\right)$,
  $\left(\begin{array}{cccccc} 0& -1& 0& 0& 1& 0 \\ -1& -1& 0& 0& 1& 0 \\ 0& -1& 0& -1& 2& 0\\-1& -2& 1& -1& 2& 0
   \\ 0& -2& 0& 0& 1& 0 \\ 0& -1& 0& 0& 0& 1
   \end{array}\right)$,\\
  $\left(\begin{array}{cccccc} 0& -1& 0& 1& -1& 1 \\ 0& -1& -1& 2& -1& 0 \\ 0& -2& -1& 2& -1& 1\\1& -2& -2& 3& -
  2& 1 \\ 1& -1& -2& 2& -1& 1 \\ 0& 0& -1& 1& -1& 1 \end{array}\right)$,
  $\left(\begin{array}{cccccc} 0& -1& 0& 1& -1& 1 \\ 0& -1& -1& 2& -1& 0 \\ 0& -1& -1& 2& -1& 1\\0& -2& -2& 3& -1&
   1 \\ 1& -1& -2& 2& -1& 1 \\ 0& 0& -1& 1& -1& 1 \end{array}\right)$,
  $\left(\begin{array}{cccccc} 0& -1& 0& 1& -1& 1 \\ 0& -1& -1& 2& -1& 0 \\ 1& -1& -1& 2& -2& 1\\1& -1& -2& 3& -2
  & 1 \\ 1& -1& -2& 2& -1& 1 \\ 0& 0& -1& 1& -1& 1
  \end{array}\right)$,\\
  $\left(\begin{array}{cccccc} 0& 0& 0& 0& -1& 1 \\ 1& 0& 0& 0& -1& 0 \\ 0& 0& 1& 0& -2& 1 \\ 1& 0& 0& 1& -3&
1 \\ 1& 1& 0& 0& -2& 1 \\ 0& 1& 0& 0& -1& 0 \end{array}\right)$,
  $\left(\begin{array}{cccccc} 0& 0& 0& 0& -1& 1 \\ 1& 0& 0& 0& -1& 0 \\ 0& 0& 1& 0& -2& 1 \\ 1& 0& 0& 1& -3&
1 \\ 1& 1& 0& 0& -2& 1 \\ 1& 0& 0& 0& -1& 1 \end{array}\right)$,
  $\left(\begin{array}{cccccc} 0& 0& 0& 0& -1& 1 \\ 1& 0& 0& 0& -1& 0 \\ 1& 1& 0& -1& -1& 1 \\ 1& 1& 1& -1& -
2& 1 \\ 1& 1& 0& 0& -2& 1 \\ 0& 1& 0& 0& -1& 0
\end{array}\right)$,\\
  $\left(\begin{array}{cccccc} 0& 0& 0& 0& -1& 1 \\ 1& 0& 0& 0& -1& 0 \\ 1& 1& 0& -1& -1& 1 \\ 1& 1& 1& -1& -
2& 1 \\ 1& 1& 0& 0& -2& 1 \\ 1& 0& 0& 0& -1& 1
\end{array}\right)$,
  $\left(\begin{array}{cccccc} 0& 1& 0& 0& -1& 1 \\ 0& 1& 0& 0& 0& 0 \\ -1& 1& 1& 0& -1& 1 \\ -1& 1& 0& 1& -1&
1 \\ -1& 1& 0& 0& 0& 1 \\ -1& 1& 0& 0& 0& 0 \end{array}\right)$,
  $\left(\begin{array}{cccccc} 0& 1& 0& 0& -1& 1 \\ 0& 1& 0& 0& 0& 0 \\ -1& 1& 1& 0& -1& 1 \\ -1& 1& 0& 1& -1&
1 \\ -1& 1& 0& 0& 0& 1 \\ 0& 0& 0& 0& 0& 1
\end{array}\right)$,\\
  $\left(\begin{array}{cccccc} 1& -1& 0& 0& 0& 0 \\ 1& 0& 0& 0& -1& 0 \\ 1& -1& 1& 0& -1& 0 \\ 2& -1& 0& 1& -
2& 0 \\ 2& 0& 0& 0& -1& 0 \\ 1& 0& 0& 0& -1& 1
\end{array}\right)$,
  $\left(\begin{array}{cccccc} 1& -1& 0& 0& 0& 0 \\ 1& 0& 0& 0& -1& 0 \\ 2& 0& 0& -1& 0& 0 \\ 2& 0& 1& -1& -1&
0 \\ 2& 0& 0& 0& -1& 0 \\ 1& 0& 0& 0& -1& 1 \end{array}\right)$,
  $\left(\begin{array}{cccccc} -1& 0& 0& 0& 0& 1 \\ -1& -1& 0& 0& 1& 0 \\ -2& -1& 1& 0& 0& 1 \\ -2& -2& 0& 1&
0& 1 \\ -1& -1& 0& 0& 0& 1 \\ -1& 0& 0& 0& 0& 0
\end{array}\right).$

The character table of $G^{s_8}$:\\
 \noindent
\begin{tabular}{c|ccccccccccccccccccccccccccc}
  & & & & & & & & & &10 & & & & & & & & & & 20    \\\hline
$\chi_8^ {(1)}$&1&1&1&1&1&1& 1& 1& 1& 1&1&1& 1& 1&1& 1&1&1 & 1&1&1    \\

$\chi_8^ {(2)}$&1&-1& 1&-1&1&-1& 1&-1& 1&-1&1&-1& 1& 1&1&-1 &-1&-1& 1      &-1&1
\\

$\chi_8^ {(3)}$&1& 1& 1& 1&1&1&-1&-1&-1&-1&-1&-1& 1& 1&1&-1 &-1&-1& 1& 1&1    \\

$\chi_8^ {(4)}$&1&-1& 1&-1&1&-1&-1& 1&-1& 1&-1&1& 1& 1&1& 1&1&1&1&-1& 1     \\

$\chi_8^ {(5)}$&1& A&-/A& /A&1&-1&-/A&/A&-A& A&1&-1&-A&-/A&1 &/A&A&-1&-/A&/A
&-A        \\

$\chi_8^ {(6)}$&1&/A&-A& A&1&-1&-A& A&-/A&/A&1&-1&-/A& -A&1& A &/A&-1&-A& A
&-/A       \\

$\chi_8^ {(7)}$&1&-/A& -A& -A&1&1& A& A&/A&/A&-1&-1&-/A& -A&1 &A&/A&-1&-A&-A
&-/A       \\
\end{tabular}

\noindent \begin{tabular}{c|ccccccccccccccccccccccccccc}
  & 22& & & & &\\\hline
$\chi_8^ {(1)}$    & 1&1&1&1 &
1&1\\

$\chi_8^ {(2)}$
&-1&
1&-1&1&1&1\\

$\chi_8^ {(3)}$   &1
&1&1&1&1&1\\

$\chi_8^ {(4)}$      &-1 &1&
-1&1&1&1\\

$\chi_8^ {(5)}$& A&1&-1&-/A&-A&-A\\

$\chi_8^ {(6)}$&/A&1&-1&-A& -/A& -/A\\

$\chi_8^ {(7)}$ &-/A&1&1&-A&-/A&-/A\\

 \ $\chi_8^ {(8)}$       &-A&1&1&-/A&-A&-A\\

$\chi_8^ {(9)}$ &/A&1&-1&-A&-/A& -/A\\

$\chi_8^ {({10})}$ & A&1&-1&-/A&-A&-A\\

$\chi_8^ {({11})}$ &-A&1&1&-/A&-A&-A\\

$\chi_8^ {({12})}$ &-/A&1&1&-A&-/A&-/A\\

$\chi_8^ {({13})}$&2 &2&2&-1&2&-1\\

$\chi_8^ {({14})}$     &
-2&2&-2&-1&2&-1\\

$\chi_8^ {({15})}$& .&-1&.&2&2&-1\\

$\chi_8^ {({16})}$       &
.&-1&.&2&2&-1\\

$\chi_8^ {({17})}$        & /B&2&2&A& /B& /A\\

$\chi_8^ {({18})}$   &B&2&2&/A&B&A\\

$\chi_8^ {({19})}$  &-/B&2&-2&A&/B&/A\\

$\chi_8^ {({20})}$ &-B&2&-2&/A&B&A\\

$\chi_8^ {({21})}$&. &-1&.&B&/B&/A\\

$\chi_8^ {({22})}$   &. &-1&.&/B&B&A\\

$\chi_8^ {({23})}$     &.&-1&.&B&/B&/A\\

$\chi_8^ {({24})}$ &.&-1&.& /B&B&A\\

$\chi_8^ {({25})}$       & .
&-2 &.&-2&4&1\\

$\chi_8^ {({26})}$   &. &-2&.&-B&/C& -/A\\

$\chi_8^ {({27})}$   & .&-2&.&-/B&C&-A\\
\end{tabular}

\noindent \noindent where A = -E(3)
  = (1-ER(-3))/2 = -b3,
B = 2*E(3)
  = -1+ER(-3) = 2b3,
C = 4*E(3)
  = -2+2*ER(-3) = 4b3.

The generators of $G^{s_9}$ are:\\
 $\left(\begin{array}{cccccc} -1& 0& 0& 0& 1& -1
\\ -1& -1& 0& 0& 1& 0
\\ -1& -1& -1& 1& 1& -1 \\ -1& -1& -1& 0& 2& -1 \\ -1& -1& 0&
0& 1& -1 \\ -1& 0& 0& 0& 0& 0
\end{array}\right)$,
$\left(\begin{array}{cccccc} -1& 1& 1& -1& 0& 0
\\ -1& 0& 2& -1& 0& 0 \\ -2& 1& 2& -1& 0& 0 \\-3& 0& 3& -1& 0& 0 \\
-2& 0& 2& 0& -1& 0
\\ -1& 0& 1& 0& 0& -1 \end{array}\right)$,
 $\left(\begin{array}{cccccc} -1& 1& 1& -1& 0& 0 \\ -1& 0& 2& -1& 0& 0 \\ -2& 0& 2& -1& 0& 0\\ -2& 0& 3& -1& -1&
 0 \\ -2& 0& 2& 0& -1& 0 \\ -1& 0& 1& 0& 0& -1 \end{array}\right).$

 The representatives   of conjugacy classes of   $G^{s_9}$ are:\\
$\left(\begin{array}{cccccc} 1& 0& 0& 0& 0& 0 \\ 0& 1& 0& 0& 0& 0 \\ 0& 0& 1& 0& 0& 0 \\
0& 0& 0& 1& 0& 0\\0& 0& 0& 0& 1& 0 \\ 0& 0& 0& 0& 0& 1
\end{array}\right)$,
  $\left(\begin{array}{cccccc} -1& 0& 0& 0& 1& -1 \\ -1& -1& 0& 0& 1& 0 \\ -1& -1& -1& 1& 1& -1\\-1& -1&
  -1& 0& 2& -1 \\ -1& -1& 0& 0& 1& -1 \\ -1& 0& 0& 0& 0& 0 \end{array}%
\right).$

The character table of $G^{s_9}$:\\
\noindent %

\noindent \noindent where A = -E(3)$^2$
  = (1+ER(-3))/2 = 1+b3,
B = -2*E(3)
  = 1-ER(-3) = 1-i3.

The generators of $G^{s_{10}}$ are:\\
 $\left(%
\right).$

The representatives   of conjugacy classes of   $G^{s_{10}}$ are:\\
$\left(%
\right).$

The character table of $G^{s_{10}}$:\\
\noindent %

The generators of $G^{s_{11}}$ are:\\
$\left(%
\right).$

The representatives   of conjugacy classes of   $G^{s_{11}}$ are:\\
$\left(%
\right).$

The character table of $G^{s_{11}}$:\\
\noindent %

\noindent \noindent where A = E(3)
  = (-1+ER(-3))/2 = b3,
B = -2*E(3)$^2$
  = 1+ER(-3) = 1+i3.

The generators of $G^{s_{12}}$ are:\\
$\left(%
\right).$

The representatives   of conjugacy classes of   $G^{s_{12}}$ are:\\
$\left(%
\right).$

The character table of $G^{s_{12}}$:\\
 \noindent %

\noindent \noindent where A = E(3)
  = (-1+ER(-3))/2 = b3,
B = 2*E(3)$^2$
  = -1-ER(-3) = -1-i3,
C = 4*E(3)$^2$
  = -2-2*ER(-3) = -2-2i3.

The generators of $G^{s_{13}}$ are:\\
$\left(%
\right).$


The character table of $G^{s_{13}}$:\\
\noindent %

A = E(5)$^4$, B = E(5)$^3$.


The generators of $G^{s_{14}}$ are:\\
 $\left(%
\right).$


The character table of $G^{s_{14}}$:\\
\noindent %

\noindent \noindent where A = -E(5)$^2$, B = -E(5)$^4$.

The generators of $G^{s_{15}}$ are:\\
$\left(%
\right).$

The character table of $G^{s_{15}}$:\\
\noindent %

\noindent \noindent where A = -E(4)
  = -ER(-1) = -i,
B = 1+E(4)
  = 1+ER(-1) = 1+i,
C = -2*E(4)
  = -2*ER(-1) = -2i,
D = -4*E(4)
  = -4*ER(-1) = -4i.

The generators of $G^{s_{16}}$ are:\\
$\left(%
\right).$

The representatives   of conjugacy classes of   $G^{s_{16}}$ are

$\left(%
\right).$


The character table of $G^{s_{16}}$:\\
\noindent %

\noindent \noindent where A = -E(8)$^3$, B = E(4)
  = ER(-1) = i.

The generators of $G^{s_{17}}$ are:\\
 $\left(%
\right).$

The character table of $G^{s_{17}}$:\\
\noindent %

\noindent \noindent where A = E(3)$^2$
  = (-1-ER(-3))/2 = -1-b3,
B = -2*E(3)$^2$
  = 1+ER(-3) = 1+i3.

The generators of $G^{s_{18}}$ are:\\
$\left(%
\right).$

The character table of $G^{s_{18}}$:\\
\noindent %

\noindent \noindent where A = E(3)
  = (-1+ER(-3))/2 = b3,
B = -E(4)
  = -ER(-1) = -i,
C = E(12)$^{11}$.

The generators of $G^{s_{19}}$ are:\\
 $\left(%
\right).$

The character table of $G^{s_{19}}$:\\
\noindent %

\noindent \noindent where A = -E(3)
  = (1-ER(-3))/2 = -b3.

The generators of $G^{s_{20}}$ are:\\
$\left(%
\right).$

The character table of $G^{s_{20}}$:\\
\noindent %

\noindent \noindent where A = E(3)
  = (-1+ER(-3))/2 = b3,
B = -E(9)$^4$-E(9)$^7$, C = E(9)$^4$, D = E(9)$^7$.

The generators of $G^{s_{21}}$ are:\\
$\left(%
\right).$

The character table of $G^{s_{21}}$:\\
{\small

\noindent %


\noindent \noindent where A = 3*E(3)$^2$
  = (-3-3*ER(-3))/2 = -3-3b3,
B = 6*E(3)$^2$
  = -3-3*ER(-3) = -3-3i3,
C = 9*E(3)$^2$
  = (-9-9*ER(-3))/2 = -9-9b3,
D = -E(3)$^2$
  = (1+ER(-3))/2 = 1+b3,
E = 2*E(3)$^2$
  = -1-ER(-3) = -1-i3,
F = -E(3)-2*E(3)$^2$
  = (3+ER(-3))/2 = 2+b3,
G = -E(3)+E(3)$^2$
  = -ER(-3) = -i3.

The generators of $G^{s_{22}}$ are:\\
$\left(%
\right).$

The character table of $G^{s_{22}}$:\\
\noindent %

\noindent \noindent where A = E(3)$^2$
  = (-1-ER(-3))/2 = -1-b3,
B = 2*E(3)$^2$
  = -1-ER(-3) = -1-i3.

The generators of $G^{s_{23}}$ are:\\
$\left(%
\right).$

The character table of $G^{s_{23}}$:\\
\noindent %

\noindent \noindent where A = E(3)
  = (-1+ER(-3))/2 = b3,
B = E(4)
  = ER(-1) = i,
C = E(12)$^7$.

The generators of $G^{s_{24}}$ are:\\
 $\left(%
\right).$

The character table of $G^{s_{24}}$:\\
{\small
\noindent %
 }

\noindent \noindent where A = E(3)$^2$
  = (-1-ER(-3))/2 = -1-b3,
B = 2*E(3)$^2$
  = -1-ER(-3) = -1-i3,
C = 3*E(3)$^2$
  = (-3-3*ER(-3))/2 = -3-3b3.

The generators of $G^{s_{25}}$ are:\\
$\left(%
\right).$

The character table of $G^{s_{25}}$:\\
\noindent %

\noindent \noindent where A = -E(3)$^2$
  = (1+ER(-3))/2 = 1+b3,
B = -2*E(3)$^2$
  = 1+ER(-3) = 1+i3.


\section {$E_7$}

 The generators of $G$ are:\\
$\left(%
\right)$

 The representatives of conjugacy classes of $G$ are:\\
$s_{1}:=\left(%
\right).$.

The character table of $G^{s_1} =G$:\\
{\small
\noindent %
}


The generators of $G^{s_2}$ are:\\
$\left(%
\right)$


The character table of $G^{s_2}$:\\
{\small
\noindent %
 }

A = E(3)$^2$
  = (-1-ER(-3))/2 = -1-b3,
B = -E(9)$^4$-E(9)$^7$,C = E(9)$^4$,D = E(9)$^7$.

The generators of $G^{s_3}$ are:\\
$\left(%
\right)$

The representatives of conjugacy classes of   $G^{s_3}$ are:\\
$\left(%
\right)$

The character table of $G^{s_3}$:\\
\noindent %

\noindent \noindent where A = E(3)$^2$
  = (-1-ER(-3))/2 = -1-b3,
B = E(9)$^7$,C = E(9)$^5$,D = -E(9)$^4$-E(9)$^7$.

The generators of $G^{s_4}$ are:\\
 $\left(%
\right)$

The representatives of conjugacy classes of   $G^{s_4}$ are:\\
$\left(%
\right)$

The character table of $G^{s_4}$:\\
${}_{
 \!\!\!\!\!\!\!\!\!\!_{
\!\!\!\!\!\!\!_{ \!\!\!\!\!\!\!  _{  \small
 \noindent
%

A = 3*E(3)$^2$
  = (-3-3*ER(-3))/2 = -3-3b3,
B = 6*E(3)$^2$
  = -3-3*ER(-3) = -3-3i3,
C = 9*E(3)$^2$
  = (-9-9*ER(-3))/2 = -9-9b3,
D = -E(3)$^2$
  = (1+ER(-3))/2 = 1+b3,
E = 2*E(3)$^2$
  = -1-ER(-3) = -1-i3,
F = E(3)+2*E(3)$^2$
  = (-3-ER(-3))/2 = -2-b3,
G = E(3)-E(3)$^2$
  = ER(-3) = i3.

The generators of $G^{s_5}$ are:\\
$\left(%
\right)$

The representatives of conjugacy classes of   $G^{s_5}$ are:\\
$\left(%
\right)$

The character table of $G^{s_5}$:\\
${}_{
 \!\!\!\!\!\!\!\!\!\!\!_{
\!\!\!\!\!\!\!_{ \!\!\!\!\!\!\!  _{  \small
 \noindent
%

A = 3*E(3)$^2$
  = (-3-3*ER(-3))/2 = -3-3b3,
B = 6*E(3)$^2$
  = -3-3*ER(-3) = -3-3i3,
C = 9*E(3)$^2$
  = (-9-9*ER(-3))/2 = -9-9b3,
D = -E(3)$^2$
  = (1+ER(-3))/2 = 1+b3,
E = 2*E(3)$^2$
  = -1-ER(-3) = -1-i3,
F = E(3)-E(3)$^2$
  = ER(-3) = i3,
G = -2*E(3)-E(3)$^2$
  = (3-ER(-3))/2 = 1-b3.

The generators of $G^{s_6}$ are:\\
$\left(%
\right)$

The character table of $G^{s_6}$:\\
\noindent %

The generators of $G^{s_7}$ are:\\
$\left(%
\right)$

The character table of $G^{s_7}$:\\
${}_{
 \!\!\!\!\!\!\!\!\!\!_{
\!\!\!\!\!\!\!_{ \!\!\!\!\!\!\!  _{  \small
 \noindent
%

A = -E(3)$^2$
  = (1+ER(-3))/2 = 1+b3,
B = 2*E(3)$^2$
  = -1-ER(-3) = -1-i3.

The generators of $G^{s_8}$ are:\\
$\left(%
\right)$

The character table of $G^{s_8}$:\\
${}_{
 \!\!\!\!\!\!_{
\!\!\!\!\!\!_{ \!\!\!\!\!\! _{  \small
 \noindent
%

A = -E(3)$^2$
  = (1+ER(-3))/2 = 1+b3,
B = 2*E(3)$^2$
  = -1-ER(-3) = -1-i3,
C = 4*E(3)$^2$
  = -2-2*ER(-3) = -2-2i3.

The generators of $G^{s_9}$ are:\\
$\left(%
\right)$

The character table of $G^{s_9}$:\\
\noindent %

The generators of $G^{s_{10}}$ are:\\
$\left(%
\right)$

The representatives of conjugacy classes of   $G^{s_{10}}$ are:\\
$\left(%
\right)$

The character table of $G^{s_{10}}$:\\
${}_{
 \!\!\!\!\!\!\!\!_{
\!\!\!\!\!\!_{ \!\!\!\!\!\!  _{  \small
 \noindent
%

\noindent where A = E(3)$^2$
  = (-1-ER(-3))/2 = -1-b3,
B = 2*E(3)$^2$
  = -1-ER(-3) = -1-i3.

The generators of $G^{s_{11}}$ are:\\
$\left(%
\right)$

The representatives of conjugacy classes of   $G^{s_{11}}$ are:\\
$\left(%
\right)$

The character table of $G^{s_{11}}$:\\
${}_{
 \!\!\!\!\!\!\!_{
\!\!\!\!\!\!\!_{ \!\!\!\!\!\!\!  _{  \small
 \noindent
%

The generators of $G^{s_{12}}$ are:\\
$\left(%
\right)$

The representatives of conjugacy classes of   $G^{s_{12}}$ are:\\
$\left(%
\right)$

The character table of $G^{s_{12}}$:\\
${}_{
 \!\!\!\!\!\!\!_{
\!\!\!\!\!\!\!_{ \!\!\!\!\!\!\!  _{  \small
 \noindent
%

A = E(3)
  = (-1+ER(-3))/2 = b3,
B = 2*E(3)
  = -1+ER(-3) = 2b3.

The generators of $G^{s_{13}}$ are:\\
$\left(%
\right)$

The character table of $G^{s_{13}}$:\\
${}_{
 \!\!\!\!\!_{
\!\!\!\!\!_{ \!\!\!\!\!  _{  \small
 \noindent
%

The generators of $G^{s_{14}}$ are:\\
$\left(%
\right)$

The character table of $G^{s_{14}}$:\\
${}_{
 \!\!\!\!\!\!\!\!\!\!\!_{
\!\!\!\!\!\!\!_{ \!\!\!\!\!\!\!  _{  \small
 \noindent
%

\noindent where A = -E(3)
  = (1-ER(-3))/2 = -b3,
B = 2*E(3)
  = -1+ER(-3) = 2b3,
C = 3*E(3)$^2$
  = (-3-3*ER(-3))/2 = -3-3b3.

The generators of $G^{s_{15}}$ are:\\
$\left(%
\right)$

The character table of $G^{s_{15}}$:\\
${}_{
 \!\!\!\!\!\!\!\!\!\!\!_{
\!\!\!\!\!\!\!_{ \!\!\!\!\!\!\!  _{  \small
 \noindent
%

\noindent \noindent where  A = E(3)$^2$
  = (-1-ER(-3))/2 = -1-b3,
B = 5*E(3)$^2$
  = (-5-5*ER(-3))/2 = -5-5b3,
C = 9*E(3)$^2$
  = (-9-9*ER(-3))/2 = -9-9b3,
D = 10*E(3)$^2$
  = -5-5*ER(-3) = -5-5i3,
E = 16*E(3)$^2$
  = -8-8*ER(-3) = -8-8i3,
F = -3*E(3)$^2$
  = (3+3*ER(-3))/2 = 3+3b3,
G = 2*E(3)$^2$
  = -1-ER(-3) = -1-i3.

The generators of $G^{s_{16}}$ are:\\
$\left(%
\right)$

The character table of $G^{s_{16}}$:\\
${}_{
 \!\!\!\!\!\!\!\!\!\!_{
\!\!\!\!\!\!\!_{ \!\!\!\!\!\!\!  _{  \small
 \noindent
%

\noindent \noindent where A = E(4)
  = ER(-1) = i,
B = 2*E(4)
  = 2*ER(-1) = 2i,
C = 3*E(4)
  = 3*ER(-1) = 3i.

The generators of $G^{s_{17}}$ are:\\
$\left(%
\right)$

The character table of $G^{s_{17}}$:\\
\noindent %

The generators of $G^{s_{18}}$ are:\\
$\left(%
\right)$

The character table of $G^{s_{18}}$:\\
${}_{
 \!\!\!\!\!\!\!_{
\!\!\!\!\!\!\!_{ \!\!\!\!\!\!\!  _{  \small
 \noindent
%

\noindent \noindent where A = E(4)
  = ER(-1) = i,
B = -1-E(4)
  = -1-ER(-1) = -1-i,
C = 2*E(4)
  = 2*ER(-1) = 2i,
D = -4*E(4)
  = -4*ER(-1) = -4i.

The generators of $G^{s_{19}}$ are:\\
$\left(%
\right)$

The character table of $G^{s_{19}}$:\\
${}_{
 \!\!\!\!\!\!\!_{
\!\!\!\!\!\!\!_{ \!\!\!\!\!\!\!  _{  \small
 \noindent
%

The generators of $G^{s_{20}}$ are:\\
$\left(%
\right)$

The character table of $G^{s_{20}}$:\\
${}_{
 \!\!\!\!\!\!\!\!\!\!\!\!\!\!\!\!_{
\!\!\!\!\!\!\!_{ \!\!\!\!\!\!\!  _{  \small
 \noindent
%

\noindent \noindent where A = E(8)$^3$, B = -E(4)
  = -ER(-1) = -i.

The generators of $G^{s_{21}}$ are:\\
$\left(%
\right)$

The character table of $G^{s_{21}}$:\\
${}_{
 \!\!\!\!\!\!\!\!\!\!\!\!_{
\!\!\!\!\!\!\!_{ \!\!\!\!\!\!\!  _{  \small
 \noindent
%

\noindent \noindent where A = -2*E(4)
  = -2*ER(-1) = -2i,
B = -4*E(4)
  = -4*ER(-1) = -4i,
C = -6*E(4)
  = -6*ER(-1) = -6i,
D = E(4)
  = ER(-1) = i,
E = -1-E(4)
  = -1-ER(-1) = -1-i.

The generators of $G^{s_{22}}$ are:\\
$\left(%
\right)$


The character table of $G^{s_{22}}$:\\
${}_{
 \!\!\!\!\!\!\!\!\!\!\!_{
\!\!\!\!\!\!\!_{ \!\!\!\!\!\!\!  _{  \small
 \noindent
%
 }}}}$

\noindent \noindent where A = E(4)
  = ER(-1) = i,
B = E(12)$^7$,
C = E(3)
  = (-1+ER(-3))/2 = b3.

The generators of $G^{s_{23}}$ are:\\
$\left(%
\right)$

The character table of $G^{s_{23}}$:\\
${}_{
 \!\!\!\!\!\!\!\!\!\!\!_{
\!\!\!\!\!\!\!_{ \!\!\!\!\!\!\!  _{  \small
\noindent
%

\noindent \noindent where A = -E(3)
  = (1-ER(-3))/2 = -b3,
B = -2*E(3)$^2$
  = 1+ER(-3) = 1+i3.

The generators of $G^{s_{24}}$ are:\\
$\left(%
\right)$

The character table of $G^{s_{24}}$:\\
${}_{
 \!\!\!\!\!\!\!\!\!\!\!\!_{
\!\!\!\!\!\!\!_{ \!\!\!\!\!\!\!  _{  \small
 \noindent
%

\noindent \noindent where A = E(4)
  = ER(-1) = i,
B = 2*E(4)
  = 2*ER(-1) = 2i,
C = 3*E(4)
  = 3*ER(-1) = 3i.

The generators of $G^{s_{25}}$ are:\\
$\left(%
\right)$

The character table of $G^{s_{25}}$:\\
${}_{
 \!\!\!\!\!\!\!\!\!\!\!\!_{
\!\!\!\!\!\!\!_{ \!\!\!\!\!\!\!  _{  \small
 \noindent
%

\noindent \noindent where A = E(4)
  = ER(-1) = i,
B = E(8).

The generators of $G^{s_{26}}$ are:\\
$\left(%
\right)$

The character table of $G^{s_{26}}$:\\
${}_{
 \!\!\!\!\!\!\!_{
\!\!\!\!\!\!\!_{ \!\!\!\!\!\!\!  _{  \small
 \noindent
%

The generators of $G^{s_{27}}$ are:\\
$\left(%
\right)$

The character table of $G^{s_{27}}$:\\
${}_{
 \!\!\!\!\!\!\!\!\!\!\!_{
\!\!\!\!\!\!\!_{ \!\!\!\!\!\!\!  _{  \small
\noindent
%

\noindent \noindent where A = E(3)$^2$
  = (-1-ER(-3))/2 = -1-b3,
B = 2*E(3)$^2$
  = -1-ER(-3) = -1-i3,
C = -3*E(3)$^2$
  = (3+3*ER(-3))/2 = 3+3b3.

The generators of $G^{s_{28}}$ are:\\
$\left(%
\right)$

The character table of $G^{s_{28}}$:\\
${}_{
 \!\!\!\!\!\!\!\!\!\!\!\!\!\!_{
\!\!\!\!\!\!\!_{ \!\!\!\!\!\!\!  _{  \small
 \noindent
%

\noindent \noindent where A = -E(4)
  = -ER(-1) = -i,
B = 2*E(4)
  = 2*ER(-1) = 2i,
C = 3*E(4)
  = 3*ER(-1) = 3i.

The generators of $G^{s_{29}}$ are:\\
$\left(%
\right)$

The character table of $G^{s_{29}}$:\\
${}_{
 \!\!\!\!\!\!\!_{
\!\!\!\!\!\!\!_{ \!\!\!\!\!\!\!  _{  \small
 \noindent
%

The generators of $G^{s_{30}}$ are:\\
$\left(%
\right)$

The character table of $G^{s_{30}}$:\\
${}_{
 \!\!\!\!\!\!\!_{
\!\!\!\!\!\!\!_{ \!\!\!\!\!\!\!  _{  \small
 \noindent
%

The generators of $G^{s_{31}}$ are:\\
$\left(%
\right)$

The character table of $G^{s_{31}}$:\\
${}_{
 \!\!\!\!\!\!\!\!_{
\!\!\!\!\!\!\!_{ \!\!\!\!\!\!\!  _{  \small
 \noindent
%

\noindent \noindent where A = -E(4)
  = -ER(-1) = -i,
B = -1+E(4)
  = -1+ER(-1) = -1+i,
C = 2*E(4)
  = 2*ER(-1) = 2i,
D = 4*E(4)
  = 4*ER(-1) = 4i.

The generators of $G^{s_{32}}$ are:\\
$\left(%
\right)$

The character table of $G^{s_{32}}$:\\
${}_{
 \!\!\!\!\!\!\!\!\!\!\!_{
\!\!\!\!\!\!\!_{ \!\!\!\!\!\!\!  _{  \small
 \noindent
%

\noindent \noindent where A = -E(3)
  = (1-ER(-3))/2 = -b3,
B = 2*E(3)
  = -1+ER(-3) = 2b3,
C = 3*E(3)
  = (-3+3*ER(-3))/2 = 3b3.

The generators of $G^{s_{33}}$ are:\\
 $\left(%
\right)$

The character table of $G^{s_{33}}$:\\
${}_{
 \!\!\!\!\!\!\!\!\!_{
\!\!\!\!\!\!\!_{ \!\!\!\!\!\!\!  _{  \small
 \noindent
%

\noindent \noindent where A = -E(3)$^2$
  = (1+ER(-3))/2 = 1+b3,
B = -2*E(3)$^2$
  = 1+ER(-3) = 1+i3.

The generators of $G^{s_{34}}$ are:\\
$\left(%
\right)$

The character table of $G^{s_{34}}$:\\
${}_{
 \!\!\!\!\!\!\!\!\!\!_{
\!\!\!\!\!\!\!_{ \!\!\!\!\!\!\!  _{  \small
 \noindent
%

\noindent \noindent where A = E(4)
  = ER(-1) = i,
B = 2*E(4)
  = 2*ER(-1) = 2i,
C = -3*E(4)
  = -3*ER(-1) = -3i.

The generators of $G^{s_{35}}$ are:\\
$\left(%
\right)$

The character table of $G^{s_{35}}$:\\
${}_{
 \!\!\!\!\!\!\!\!\!\!\!_{
\!\!\!\!\!\!\!_{ \!\!\!\!\!\!\!  _{  \small
 \noindent
%
 }}}}$

\noindent \noindent where A = -E(3)
  = (1-ER(-3))/2 = -b3,
B = -E(4)
  = -ER(-1) = -i,
C = -E(12)$^7$.

The generators of $G^{s_{36}}$ are:\\
$\left(%
\right)$

The character table of $G^{s_{36}}$:\\
${}_{
 \!\!\!\!\!\!\!\!\!\!\!_{
\!\!\!\!\!\!\!_{ \!\!\!\!\!\!\!  _{
 \noindent
%

\noindent \noindent where A = E(3)$^2$
  = (-1-ER(-3))/2 = -1-b3,
B = 2*E(3)$^2$
  = -1-ER(-3) = -1-i3.

The generators of $G^{s_{37}}$ are:\\
$\left(%
\right)$

The character table of $G^{s_{37}}$:\\
${}_{
 \!\!\!\!\!\!\!_{
\!\!\!\!\!\!\!_{ \!\!\!\!\!\!\!  _{  \small
 \noindent
%

\noindent \noindent where A = E(3)
  = (-1+ER(-3))/2 = b3,
B = 2*E(3)$^2$
  = -1-ER(-3) = -1-i3,
C = -3*E(3)$^2$
  = (3+3*ER(-3))/2 = 3+3b3.

The generators of $G^{s_{38}}$ are:\\
$\left(%
\right)$

The character table of $G^{s_{38}}$:\\
${}_{
 \!\!\!\!\!\!\!\!\!_{
\!\!\!\!\!\!\!_{ \!\!\!\!\!\!\!  _{  \small
 \noindent
%
\noindent \noindent where A = E(3)
  = (-1+ER(-3))/2 = b3,
B = 2*E(3)
  = -1+ER(-3) = 2b3,
C = 4*E(3)$^2$
  = -2-2*ER(-3) = -2-2i3.

The generators of $G^{s_{39}}$ are:\\
$\left(%
\right)$

The character table of $G^{s_{39}}$:\\
${}_{
 \!\!\!\!\!\!\!_{
\!\!\!\!\!\!\!_{ \!\!\!\!\!\!\!  _{  \small
 \noindent
%
 }}}}$

 \noindent \noindent where A = E(3)$^2$
  = (-1-ER(-3))/2 = -1-b3.

The generators of $G^{s_{40}}$ are:\\
$\left(%
\right)$

The character table of $G^{s_{40}}$:\\
\noindent %

\noindent \noindent where A = -E(5),B = -E(5)$^2$.

The generators of $G^{s_{41}}$ are:\\
$\left(%
\right)$

The character table of $G^{s_{41}}$:\\
${}_{
 \!\!\!\!\!\!\!\!\!\!\!\!\!\!\!\!_{
\!\!\!\!\!\!\!_{ \!\!\!\!\!\!\!  _{  \small
 \noindent
%

\noindent \noindent where A = E(5),B = E(5)$^2$,C = 2*E(5)$^4$,D =
2*E(5)$^3$.

The generators of $G^{s_{42}}$ are:\\
$\left(%
\right)$

The character table of $G^{s_{42}}$:\\
${}_{
 \!\!\!\!\!\!\!\!\!\!\!_{
\!\!\!\!\!\!\!_{ \!\!\!\!\!\!\!  _{  \small
 \noindent
%
 }}}}$

\noindent \noindent where A = -E(3)
  = (1-ER(-3))/2 = -b3,
B = -E(4)
  = -ER(-1) = -i,
C = -E(12)$^7$.

The generators of $G^{s_{43}}$ are:\\
$\left(%
\right)$

The character table of $G^{s_{43}}$:\\
${}_{
 \!\!\!\!\!\!\!_{
\!\!\!\!\!\!\!_{ \!\!\!\!\!\!\!  _{  \small
 \noindent
%

\noindent \noindent where A = -E(3)$^2$
  = (1+ER(-3))/2 = 1+b3,
B = 2*E(3)
  = -1+ER(-3) = 2b3,
C = 3*E(3)
  = (-3+3*ER(-3))/2 = 3b3.

The generators of $G^{s_{44}}$ are:\\
$\left(%
\right)$

The character table of $G^{s_{44}}$:\\
${}_{
 \!\!\!\!\!\!\!\!\!\!\!\!\!\!\!\!\!_{
\!\!\!\!\!\!\!_{ \!\!\!\!\!\!\!  _{  \small
 \noindent
%

\noindent \noindent where A = E(4)
  = ER(-1) = i.

The generators of $G^{s_{45}}$ are:\\
$\left(%
\right)$

The character table of $G^{s_{45}}$:\\
${}_{
 \!\!\!\!\!\!\!\!\!\!\!_{
\!\!\!\!\!\!\!_{ \!\!\!\!\!\!\!  _{  \small
 \noindent
%

\noindent \noindent where A = E(3)$^2$
  = (-1-ER(-3))/2 = -1-b3,
B = 5*E(3)$^2$
  = (-5-5*ER(-3))/2 = -5-5b3,
C = 9*E(3)$^2$
  = (-9-9*ER(-3))/2 = -9-9b3,
D = 10*E(3)$^2$
  = -5-5*ER(-3) = -5-5i3,
E = 16*E(3)$^2$
  = -8-8*ER(-3) = -8-8i3,
F = -3*E(3)$^2$
  = (3+3*ER(-3))/2 = 3+3b3,
G = 2*E(3)$^2$
  = -1-ER(-3) = -1-i3.

The generators of $G^{s_{46}}$ are:\\
$\left(%
\right)$

The character table of $G^{s_{46}}$:\\
${}_{
 \!\!\!\!\!\!\!\!\!_{
\!\!\!\!\!\!\!_{ \!\!\!\!\!\!\!  _{  \small
 \noindent
%
 }}}}$

\noindent \noindent where A = -E(3)
  = (1-ER(-3))/2 = -b3.

The generators of $G^{s_{47}}$ are:\\
$\left(%
\right)$

The character table of $G^{s_{47}}$:\\
\noindent %

\noindent \noindent where A = -E(5),B = -E(5)$^2$.

The generators of $G^{s_{48}}$ are:\\
$\left(%
\right)$

The character table of $G^{s_{48}}$:\\
${}_{
 \!\!\!\!\!\!\!\!\!\!\!_{
\!\!\!\!\!\!\!_{ \!\!\!\!\!\!\!  _{  \small
 \noindent
%

A = E(5),B = E(5)$^2$,C = -2*E(5)$^4$,D = -2*E(5)$^3$.

The generators of $G^{s_{49}}$ are:\\
$\left(%
\right)$,

The character table of $G^{s_{49}}$:\\
${}_{
 \!\!\!\!\!\!\!\!\!\!\!\!_{
\!\!\!\!\!\!\!_{ \!\!\!\!\!\!\!  _{  \small
 \noindent
%

\noindent \noindent where A = E(3)$^2$
  = (-1-ER(-3))/2 = -1-b3,
B = E(5),C = E(15)$^8$,D = E(15)$^{13}$,E = E(5)$^2$, F
=E(15)$^{11}$,G = E(15).

The generators of $G^{s_{50}}$ are:\\
$\left(%
\right)$

The character table of $G^{s_{50}}$:\\
${}_{
 \!\!\!\!\!\!\!_{
\!\!\!\!\!\!\!_{ \!\!\!\!\!\!\!  _{  \small
 \noindent
%

\noindent \noindent where A = -E(3)
  = (1-ER(-3))/2 = -b3,
B = -E(5),C = -E(5)$^2$ ,D = -E(15)$^{13}$,E = -E(15)$^8$,F =
-E(15)$^{11}$,G = -E(15).

The generators of $G^{s_{51}}$ are:\\
$\left(%
\right)$

The character table of $G^{s_{51}}$:\\

\noindent %

\noindent \noindent where A = -E(7)$^6$,B = -E(7)$^5$,C = -E(7)$^4$.

The generators of $G^{s_{52}}$ are:\\
$\left(%
\right)$

The character table of $G^{s_{52}}$:\\
\noindent %

\noindent \noindent where A = E(7)$^4$,B = E(7),C = E(7)$^5$.

The generators of $G^{s_{53}}$ are:\\
$\left(%
\right)$

The character table of $G^{s_{53}}$:\\
${}_{
 \!\!\!\!\!\!\!\!\!_{
\!\!\!\!\!\!\!_{ \!\!\!\!\!\!\!  _{  \small
 \noindent
%

\noindent \noindent where A = -E(3)
  = (1-ER(-3))/2 = -b3,
B = 2*E(3)$^2$
  = -1-ER(-3) = -1-i3,
C = -3*E(3)$^2$
  = (3+3*ER(-3))/2 = 3+3b3.

The generators of $G^{s_{54}}$ are:\\
$\left(%
\right)$

The character table of $G^{s_{54}}$:\\
${}_{
 \!\!\!\!\!\!\!\!\!_{
\!\!\!\!\!\!\!_{ \!\!\!\!\!\!\!  _{  \small
 \noindent
%

\noindent \noindent where A = E(4)
  = ER(-1) = i,
B = -E(8).

The generators of $G^{s_{55}}$ are:\\
$\left(%
\right)$

The character table of $G^{s_{55}}$:\\
${}_{
 \!\!\!\!\!\!\!\!\!_{
\!\!\!\!\!\!\!_{ \!\!\!\!\!\!\!  _{  \small
 \noindent
%
 }}}}$

\noindent \noindent where A = -E(4)
  = -ER(-1) = -i,
B = -E(12)$^7$,C = E(3)
  = (-1+ER(-3))/2 = b3.

The generators of $G^{s_{56}}$ are:\\
$\left(%
\right)$

The character table of $G^{s_{56}}$:\\
${}_{
 \!\!\!\!\!\!\!_{
\!\!\!\!\!\!\!_{ \!\!\!\!\!\!\!  _{  \small
 \noindent
%

\noindent \noindent where A = E(4)
  = ER(-1) = i,
B = E(8).

The generators of $G^{s_{57}}$ are:\\
$\left(%
\right)$

The character table of $G^{s_{57}}$:\\
${}_{
 \!\!\!\!\!\!\!\!\!\!_{
\!\!\!\!\!\!\!_{ \!\!\!\!\!\!\!  _{  \small
 \noindent
%

\noindent \noindent where A = -2*E(4)
  = -2*ER(-1) = -2i,
B = -4*E(4)
  = -4*ER(-1) = -4i,
C = -6*E(4)
  = -6*ER(-1) = -6i,
D = E(4)
  = ER(-1) = i,
E = -1+E(4)
  = -1+ER(-1) = -1+i.

The generators of $G^{s_{58}}$ are:\\
$\left(%
\right)$

The character table of $G^{s_{58}}$:\\
${}_{
 \!\!\!\!\!\!\!_{
\!\!\!\!\!\!\!_{ \!\!\!\!\!\!\!  _{  \small
 \noindent
%
 }}}}$

\noindent \noindent where A = -E(4)
  = -ER(-1) = -i,
B = -E(12)$^7$,
C = E(3)
  = (-1+ER(-3))/2 = b3.

The generators of $G^{s_{59}}$ are:\\
$\left(%
\right)$

The character table of $G^{s_{59}}$:\\
${}_{
 \!\!\!\!\!\!\!_{
\!\!\!\!\!\!_{ \!\!\!\!\!\!  _{  \small
 \noindent
%
 }}}}$

\noindent \noindent where A = -E(4)
  = -ER(-1) = -i,
B = E(3)$^2$
  = (-1-ER(-3))/2 = -1-b3,
C = E(12)$^{11}$.

The generators of $G^{s_{60}}$ are:\\
$\left(%
\right)$

The character table of $G^{s_{60}}$:\\
${}_{
 \!\!\!\!\!\!\!\!\!\!\!\!\!\!\!\!\!\!_{
\!\!\!\!\!\!\!_{ \!\!\!\!\!\!\!  _{  \small
 \noindent
%

\noindent \noindent where A = E(4)
  = ER(-1) = i.

\section {$F_4$}

The generators of $G$ are:\\
$\left(\begin{array}{cccc}-1&0&1&0 \\
0&1&0&0 \\ 0&0&1&0
\\ 0&0&0&1 \end{array}%
\right)$,

The representatives of conjugacy classes of $G$ are:\\
$s_{1}:=\left(%
\right)$

The character table of $G^{s_1} =G$:\\
\noindent %

The generators of $G^{s_2}$ are:\\
$\left(%
\right)$

The character table of $G^{s_2}$:\\
\noindent %

The generators of $G^{s_3}$ are:\\
$\left(%
\right)$

The character table of $G^{s_3}$:\\
\noindent %

The generators of $G^{s_4}$ are:\\
$\left(%
\right)$

The character table of $G^{s_4}$:\\
\noindent %

\noindent \noindent where A = E(4)
  = ER(-1) = i,
B = -1+E(4)
  = -1+ER(-1) = -1+i,
C = 2*E(4)
  = 2*ER(-1) = 2i,
D = 4*E(4)
  = 4*ER(-1) = 4i.

The generators of $G^{s_5}$ are:\\
$\left(%
\right)$

The character table of $G^{s_5}$:\\
\noindent %

\noindent \noindent where A = E(3)$^2$
  = (-1-ER(-3))/2 = -1-b3,
B = 2*E(3)$^2$
  = -1-ER(-3) = -1-i3.

The generators of $G^{s_6}$ are:\\
$\left(%
\right)$

The character table of $G^{s_6}$:\\
\noindent %

\noindent \noindent where A = E(3)$^2$
  = (-1-ER(-3))/2 = -1-b3,
B = 2*E(3)$^2$
  = -1-ER(-3) = -1-i3.

The generators of $G^{s_7}$ are:\\
$\left(%
\right)$

The character table of $G^{s_7}$:\\
\noindent %

The generators of $G^{s_8}$ are:\\
$\left(%
\right)$

The character table of $G^{s_8}$:\\
\noindent %

The generators of $G^{s_9}$ are:\\
$\left(%
\right)$)

 The representatives   of conjugacy classes of   $G^{s_9}$ are:\\

$\left(%
\right)$

The character table of $G^{s_9}$:\\
\noindent %

\noindent \noindent where A = E(4)
  = ER(-1) = i.

The generators of $G^{s_{10}}$ are:\\
$\left(%
\right)$

The character table of $G^{s_{10}}$:\\
\noindent %

\noindent \noindent where A = -E(3)$^2$
  = (1+ER(-3))/2 = 1+b3,
B = -2*E(3)
  = 1-ER(-3) = 1-i3.

The generators of $G^{s_{11}}$ are:\\
$\left(%
\right)$

The character table of $G^{s_{11}}$:\\
\noindent %

\noindent \noindent where A = -E(3)$^2$
  = (1+ER(-3))/2 = 1+b3,
B = -2*E(3)
  = 1-ER(-3) = 1-i3.

The generators of $G^{s_{12}}$ are:\\
$\left(%
\right)$

The character table of $G^{s_{12}}$:\\
\noindent %

\noindent \noindent where A = E(3)
  = (-1+ER(-3))/2 = b3,
B = 2*E(3)$^2$
  = -1-ER(-3) = -1-i3,
C = 3*E(3)$^2$
  = (-3-3*ER(-3))/2 = -3-3b3.

The generators of $G^{s_{13}}$ are:\\
$\left(%
\right)$

The character table of $G^{s_{13}}$:\\
\noindent %

\noindent \noindent where A = E(3)
  = (-1+ER(-3))/2 = b3,
B = 2*E(3)$^2$
  = -1-ER(-3) = -1-i3,
C = 3*E(3)$^2$
  = (-3-3*ER(-3))/2 = -3-3b3.

The generators of $G^{s_{14}}$ are:\\
$\left(%
\right)$

The character table of $G^{s_{14}}$:\\
\noindent %

\noindent \noindent where A = -E(3)$^2$
  = (1+ER(-3))/2 = 1+b3,
B = -E(4)
  = -ER(-1) = -i,
C = -E(12)$^{11}$.

The generators of $G^{s_{15}}$ are:\\
$\left(%
\right)$

The character table of $G^{s_{15}}$:\\
\noindent %

\noindent \noindent where A = -E(3)
  = (1-ER(-3))/2 = -b3.

The generators of $G^{s_{16}}$ are:\\
$\left(%
\right)$

The character table of $G^{s_{16}}$:\\
\noindent %

\noindent \noindent where A = -E(3)
  = (1-ER(-3))/2 = -b3.

The generators of $G^{s_{17}}$ are:\\
$\left(%
\right)$

The character table of $G^{s_{17}}$:\\
\noindent %

The generators of $G^{s_{18}}$ are:\\
$\left(%
\right)$

The character table of $G^{s_{18}}$:\\
\noindent %

The generators of $G^{s_{19}}$ are:\\
$\left(%
\right)$

The character table of $G^{s_{19}}$:\\
\noindent %

\noindent \noindent where A = E(4)
  = ER(-1) = i.

The generators of $G^{s_{20}}$ are:\\
$\left(%
\right)$

The character table of $G^{s_{20}}$:\\
\noindent %

\noindent \noindent where A = E(3)$^2$
  = (-1-ER(-3))/2 = -1-b3.

The generators of $G^{s_{21}}$ are:\\
$\left(%
\right)$

The character table of $G^{s_{21}}$:\\
\noindent %

\noindent \noindent where A = E(3)$^2$
  = (-1-ER(-3))/2 = -1-b3.

The generators of $G^{s_{22}}$ are:\\
$\left(%
\right)$

The character table of $G^{s_{22}}$:\\
\noindent %

The generators of $G^{s_{23}}$ are:\\
$\left(%
\right)$

The character table of $G^{s_{23}}$:\\
\noindent %

\noindent \noindent where A = E(8)$^3$, B = E(4)
  = ER(-1) = i.

The generators of $G^{s_{24}}$ are:\\
$\left(%
\right)$

The character table of $G^{s_{24}}$:\\
\noindent %

\noindent \noindent where A = E(4)
  = ER(-1) = i,
B = 2*E(4)
  = 2*ER(-1) = 2i.

The generators of $G^{s_{25}}$ are:\\
$\left(%
\right)$

The character table of $G^{s_{25}}$:\\
\noindent %

\noindent \noindent where A = E(4)
  = ER(-1) = i,
B = 2*E(4)
  = 2*ER(-1) = 2i.

\section {$G_2$}

The generators of $G$ are:\\
$\left(\begin{array}{cc} -1&1 \\ 0&1
\end{array}%
\right)$

The representatives of conjugacy classes of $G$ are:\\
$s_{1}:=\left(%
\right).$

The character table of $G^{s_1} =G$:\\
\noindent %

The generators of $G^{s_2}$ are:\\
$\left(%
\right)$

The character table of $G^{s_2}$:\\
\noindent %

The generators of $G^{s_3}$ are:\\
$\left(%
\right)$

The character table of $G^{s_3}$:\\
\noindent %

The generators of $G^{s_4}$ are:\\
$\left(%
\right)$

The character table of $G^{s_4}$:\\
\noindent %

The generators of $G^{s_5}$ are:\\
$\left(%
\right)$

The character table of $G^{s_5}$:\\
\noindent %

\noindent \noindent where A = -E(3)$^2$
  = (1+ER(-3))/2 = 1+b3.

The generators of $G^{s_6}$ are:\\
$\left(%
\right)$

The character table of $G^{s_6}$:\\
\noindent %

\noindent \noindent where A = -E(3)$^2$
  = (1+ER(-3))/2 = 1+b3.

\vskip 0.3cm

\vskip 0.3cm

\noindent{\large\bf Acknowledgement}: We would like to thank Prof.
N. Andruskiewitsch and Dr. F. Fantino for suggestions and help.

\end {document}